\documentclass[10pt,leqno]{amsart} 

   \usepackage[centertags]{amsmath}
   \usepackage{amsfonts}
   \usepackage{amsmath}
   \usepackage{amssymb}
   \usepackage{amsthm}
   \usepackage{newlfont}
   \usepackage[latin1]{inputenc}
\usepackage{xcolor}

\allowdisplaybreaks[3]

\theoremstyle{plain}
\newtheorem{Theo}{Theorem}[section]

\newtheorem{Cor}[Theo]{Corollary}
\newtheorem{Lem}[Theo]{Lemma}
\newtheorem{TheoA}{Theorem}

\theoremstyle{definition}
\newtheorem{Def}{Definition}[section]

\theoremstyle{remark}
\newtheorem{Rem}{Remark}[section]

\newcommand{\Ea}{E_{\mu}}
\newcommand{\La}{\Lambda_{\mu}}

\newcommand{\I}{\mathcal{I}}

\newcommand{\lp}{\mathcal{L}\text{-}\mathcal{P}}
\newcommand{\psa}{JPS }

\newcommand{\RR}{\mathbb{R}}

\newcommand{\NN}{\mathbb{N}}
\newcommand{\CC}{\mathbb{C}}

\newcommand{\rrp}{\mathcal{R}}

\newcommand{\cc}{\gamma}

\numberwithin{equation}{section}

\newcommand\sign{\operatorname{sign}}

\newcommand*\pFq{\begingroup
        \dopFq
}
\catcode`\,12
\def\dopFq#1#2#3#4#5{%
        {}_{#1}F_{#2}\biggl(\genfrac..{0pt}{}{#3}{#4};#5\biggr)%
        \endgroup
}

\title{Brenke polynomials with real zeros and the Riemann Hypothesis}
\author{Antonio J. Dur\'an}
\address{Departamento de An\'a\-li\-sis Mate\-m\'a\-ti\-co and IMUS,
        Universidad de Sevilla,
        41080 Sevilla, Spain}
\email{duran@us.es}
   \date{}

   \thanks{This research was partially supported by PID2021-124332NB-C21
(Minis\-te\-rio de Cien\-cia e Inno\-va\-ci\'on and Feder Funds (European Union)), and
FQM-262 (Jun\-ta de Anda\-lu\-c\'ia).}



\begin{document}

\begin{abstract}
If $A(z)=\sum_{n=0}^\infty a_nz^n$ and $B(z)=\sum_{n=0}^\infty b_nz^n$ are two formal power series, with $a_n,b_n\in \RR$, the polynomials $(p_n)_n$ defined by the generating function
$$
A(z)B(xz)=\sum_{n=0}^\infty p_n(x)z^n
$$
are called the Brenke polynomials generated by $A$ and associated to $B$.
We say that $A\in \rrp_B$ if the Brenke polynomials $(p_n)_n$ have only real zeros. Among other results, in this paper we find necessary and sufficient conditions on $B$ such that $\rrp_B=\lp$, where $\lp$ denotes the Laguerre-Pólya class (of entire functions). These results can be considered an extension to Brenke polynomials of the Jensen, and Pólya and Schur characterization $\rrp_{e^z}=\lp$, for Appell polynomials. When applying our results to a relative of the Riemann zeta function, we find new equivalencies for the Riemann Hypothesis in terms of real-rootedness of some sequences of Brenke polynomials.
\end{abstract}

   \maketitle

\section{Introduction and results}
Let $A$ and $B$ be two formal power series
$$
A(z)=\sum_{n=0}^\infty a_nz^n,\quad B(z)=\sum_{n=0}^\infty b_nz^n,
$$
with $a_n,b_n\in \RR$, and normalized by taking $a_0=b_0=1$.

We say that the polynomials $(p_n)_n$ are Brenke polynomials generated by $A$ and associated to $B$ if
\begin{equation}\label{dgf}
A(z)B(xz)=\sum_{n=0}^\infty p_n(x)z^n.
\end{equation}

The case when $B(z)=e^z$ are the well-known Appell polynomials.

In this paper we pose the (so-called) real-rooted Brenke polynomial problem: that is, the problem of determining when the Brenke polynomials generated by $A$ and associated to $B$ have only real zeros. Hence we associate to $B$ the set  $\rrp_B$ defined by
\begin{align}\label{rrppi}
\rrp _B&=\{\mbox{$A$: $A$ is a formal power series, $a_0=1$, such that }\\\nonumber
&\hspace{1.2cm}\mbox{the Brenke polynomials generated by $A$ have only real zeros}\}.
\end{align}
There is other related set $\rrp_B^p$ defined by
\begin{align}\label{rrpp2i}
\rrp _B^p&=\{\mbox{$A$: $A$ is a formal power series, $a_0=1$, and infinitely many of}\\\nonumber
&\hspace{1.2cm}\mbox{the Brenke polynomials generated by $A$ have only real zeros}\}.
\end{align}

The characterization of $\rrp_B$ for the case $B(z)=e^z$, i.e., the Appell case, goes more than one century back and it is due, independently, to Jensen and also to Pólya and Schur. This case gave rise to important results and concepts.

\begin{Def}\label{def1} An entire function $A$ is said to be in the Laguerre-Pólya class if it can be expressed in the form
\begin{equation}\label{pspr}
A(z)=cz^me^{-az^2+bz}\prod_{j=1}^\infty \left(1-\zeta_j z\right)e^{\zeta_j z},
\end{equation}
where $a\ge 0$, $m\in \NN$, $b, c,\zeta_j\in\RR$, $j\ge 1$, and $\sum_{j=1}^\infty \zeta_j^2<+\infty$.
The Laguerre-Pólya class will be denoted by $\lp$.
\end{Def}

The Laguerre-Pólya class characterizes the formal power series generating Appell polynomials with only real zeros.

\begin{TheoA}\label{TPS} Let $A$ be a formal power series with $a_0=1$. The Appell polynomials generated by $A$ have all their zeros real if and only if $A$ is an entire function in the Laguerre-Pólya class. In other words,
$\rrp_{e^z}=\{A: \mbox{$A\in\lp$ and $A(0)=1$}\}$.
\end{TheoA}

See \cite{jen} (although Jensen used some analyticity assumption on $A$), \cite{PS}, and also \cite{CrCs0}.

Actually, Jensen considered the polynomials $q_n(z)=z^np_n(1/z)$,
where $(p_n)_n$ are the Appell polynomials generated by $A$, that is
\begin{equation}\label{jepo}
q_n(z)=\sum_{j=0}^n\frac{a_j}{(n-j)!}z^j.
\end{equation}
In particular, he proved that
\begin{equation}\label{basy}
\lim_nn!\,q_n(z/n)=A(z)
\end{equation}
uniformly in compact sets of $\CC$.

Since Appell polynomials satisfy $p_n'=p_{n-1}$, it is easy to deduce from Theorem~\ref{TPS} that also $\rrp_{e^z}^p=\rrp_{e^z}=\{A: \mbox{$A\in\lp$ and $A(0)=1$}\}$.

We also consider functions of first type in the Laguerre-Pólya class.

\begin{Def}
We say that an entire function $A$ in the Laguerre-Pólya class is of type I (or first type), in short $A\in \lp I$, if
$A(z)$ or $A(-z)$ has a product representation of the form
$$
cz^me^{\alpha z}\prod_{k=1}^\infty \left(1+\zeta_kz\right),
$$
where $\alpha \ge 0$, $c\in \RR$, $m\in \NN$ and $\zeta_k> 0$, $\sum_k\zeta_k<\infty$.
\end{Def}

The class $\lp I$ has been extensively studied, mainly because its usefulness in the study of polynomials with real zeros.

There is also an important connection between the classes $\lp$ and $\lp I$ and the Riemann Hypothesis (RH in short). Indeed, let $\xi$ be the Riemann xi function
$$
\xi(s)=\frac{1}{2}s(s-1)\pi^{-s/2}\Gamma\left(\frac{s}{2}\right)\zeta(s),
$$
where $\Gamma$ is the Gamma function and $\zeta$ is the Riemann zeta function. Define then
\begin{equation}\label{var}
\varsigma(z)=\frac{1}{\xi(1/2)}\xi\left(\frac{1}{2}+\sqrt z\right)=\sum_{n=0}^\infty \frac{\gamma_n}{n!}z^n,
\end{equation}
where
\begin{equation}\label{varg}
\gamma_n=\frac{n!\,\xi^{(2n)}(1/2)}{(2n)!\,\xi(1/2)}.
\end{equation}

The Riemann Hypothesis can be formulated in terms of the zeros of $\varsigma$:  the Riemann Hypothesis is true if and only if the function $\varsigma$ has only real zeros. Since $\varsigma$ has order $1/2$, on the one hand we have that RH is equivalent to  $\varsigma\in \lp $. On the other hand, since $\gamma_n>0$ (see Corollary \ref{csi} below), then RH is also equivalent to $\varsigma\in \lp I$.
According to Theorem \ref{TPS}, and as pointed out by Jensen, if we write $(q_n)_n$ for the Jensen polynomials associated to $\varsigma$ so that
\begin{equation}\label{jenp}
q_n(z)=\sum_{j=0}^n\frac{\gamma_j}{(n-j)!\,j!}z^j,
\end{equation}
then RH is equivalent to the polynomials $q_n$, $n\ge 0$, having only real zeros.
Or, in terms of the Appell polynomials generated by $\varsigma$: the polynomials $p_n(x)=x^nq_n(1/x)$ has only real zeros for $n\ge 0$ if and only if RH is true.

The Laguerre-Pólya class has other interesting characterization due to Laguerre and Pólya: the functions in the Laguerre-Pólya class are those entire functions which can be approximated (locally uniformly) by polynomials with only real zeros.

\begin{TheoA}\label{thiwi}
An analytic function $A$ at $z=0$ is in the Laguerre-Pólya class if and only if there exist polynomials $r_n$ with only real zeros such that
$$
\lim_{n\to \infty}r_n(z)=A(z)
$$
uniformly in $\vert z\vert \le r$ for some $r>0$.
\end{TheoA}

See \cite[Theorem 3.3, p. 43]{HiWi}; the implication $\Leftarrow$ was proved by Laguerre and $\Rightarrow$ by Pólya. As an easy  consequence of Theorem \ref{thiwi} it follows that the Laguerre-Pólya class is closed under differentiation.

The main problem we study in this paper is for which formal power series $B$, the Jensen and Pólya-Schur characterization result for Appell polynomials can be extended for the Brenke polynomials associated to $B$. That is, we
find necessary and sufficient conditions on the formal power series $B$ such that $\rrp_B=\{A: \mbox{$A\in\lp$ and $A(0)=1$}\}$.
When that happens, we say that $B$ belongs to the \psa class (\textbf{J} for Jensen, \textbf{P} for Pólya and \textbf{S} for Schur).

The content of the paper is as follows.

In Section \ref{bp} we consider some basic facts and new results about Brenke polynomials.
In particular, we prove an asymptotic property for Brenke polynomials under mild conditions on $A$ and $B$ (which can be considered an analog to the asymptotic (\ref{basy}) for the Jensen polynomials). In the case of $B$, we assume that
\begin{equation}\label{pbsxxi}
\mbox{$b_n\not=0$, $n$ big enough, and } \lim_{n\to \infty}\frac{b_{n-2}b_n}{b_{n-1}^2}=1.
\end{equation}
The existence of this limit is equivalent to the existence of a sequence $(\tau_n)_n$, $\tau_n\not =0$, such that
\begin{equation}\label{pbsi}
\lim_{n\to \infty}\frac{b_{n-j}}{b_n\tau_n^j}=1,\quad \mbox{for all $j\ge 1$}.
\end{equation}
We then have the following Theorem.

\begin{Theo}\label{iasy}
Assume that the sequence $(b_n)_n$ satisfies (\ref{pbsxxi}), and let $(\tau_n)_n$ be a sequence satisfying (\ref{pbsi}).
For an analytic function $A$ at $z=0$ (with radius of convergence $r$) consider the Brenke polynomials $(p_n)_n$ generated by $A$ and associated to $B$.
If there exist $M,N >0$ such that
\begin{equation}\label{cbt}
\left|\frac{b_{n-j}}{b_n\tau_n^j}\right|\le M,\quad \mbox{for $j,n\ge N$},
\end{equation}
then
\begin{equation}\label{las2}
\lim_{n\to \infty}\left(\frac{z}{\tau_n}\right)^n\frac{p_{n}(\tau_n/z)}{b_n}=\sum_{j=0}^{\infty}a_jz^{j}=A(z)
\end{equation}
uniformly in compact sets of $\{z:\vert z\vert <r\}$.
\end{Theo}

In Section \ref{rrbp}, we pose the real-rooted Brenke polynomial problem (of determining the set $\rrp_B$ for a formal power series $B$ (\ref{rrppi})) and consider a couple of illustrative examples (one of them the Appell case). We have not found this problem explicitly posed in the literature, although
the problem has been studied and solved for some particular $B$'s (v.g., $B(z)=e^z$) and it is related to some others significant problems studied during the last 100 years (such as the problem of characterizing power series having sections with only real zeros). There are also some related results in \cite[Chapter 2]{ili} (see also the references therein).

In Section \ref{jps} we study the \psa class.

Along this paper, we use a normalized version of the classes $\lp$ and $\lp I$:
\begin{align}\label{lpn}
\lp_0&=\{A: \mbox{$A\in \lp$ and $A(0)=1$}\},\\\label{lpn0}
\lp_0 I&=\{A: \mbox{$A\in \lp I$ and $A(0)=1$}\}.
\end{align}

As one of the main results in this paper, we characterize the class $\lp I$ in terms of real-rootedness of the Brenke polynomials associated to $B$.

\begin{Theo} \label{bpoli}
Let $B$ be a formal power series with $b_0=1$. The following conditions are equivalent.
\begin{enumerate}
\item $B$ is an entire function of the first type in the Laguerre-Pólya class.
\item $\lp_0 \subset \rrp _B$.
\end{enumerate}
\end{Theo}

According to Theorem \ref{bpoli}, the Riemann Hypothesis is equivalent to $\lp_0\subset \rrp_\varsigma$. Moreover, we also have (compare with (\ref{jenp})).

\begin{Cor}\label{rhe} For a non-negative integer $N$, the following conditions are equivalent.
\begin{enumerate}
\item The Riemann Hypothesis is true.
\item The polynomials
\begin{equation}\label{jenpx}
\hat q_{N,n}(x)=\sum_{j=0}^n\frac{\gamma_j}{(n+j)!^Nj!\,(n-j)!}x^j
\end{equation}
have only real zeros for all $n\ge 0$.
\end{enumerate}
\end{Cor}
\medskip

Actually, if we take $N$ to be a non-negative real number, our computations show that the polynomials $\hat q_{N,n}$ still seem to have only real zeros.

\medskip

Theorem \ref{bpoli} establishes that  $B\in \lp I$ is a necessary condition for a formal power series $B$ to be in the \psa class.
We have also proved some sufficient conditions on a formal series $B$ so that $\rrp _B=\lp_0 $.

\begin{Theo}\label{prin2} Assume $B\in \lp I$, $B(0)=1$, it is not a polynomial and
\begin{equation}\label{supl}
\lim_{n\to \infty} \frac{b_{n-2}b_n}{b_{n-1}^2}=1.
\end{equation}
Then $\rrp_B=\lp_0$.
\end{Theo}

In order to prove Theorem \ref{prin2}, we use the asymptotic in Theorem \ref{iasy} and the following improvement of Theorem \ref{thiwi}:

\begin{Theo}\label{tlli}
Let $(n_k)_k$ be an increasing sequence of positive integers. Assume that we have polynomials $r_{n_k}(z)=\sum_{j=0}^{n_k}a_{j,n_k}z^j$, $k\ge 0$, with only real zeros, and that the following limits exist for all $j\ge 0$:
$$
\lim_ka_{j,n_k}=a_j.
$$
Then $A(z)=\sum_{j=0}^\infty a_{j}z^j$ is an entire function in the Laguerre-Pólya class.
\end{Theo}

It follows easily from the Grosswald asymptotic for $\xi^{(2n)}(1/2)$ (see \cite{gro1}, \cite{gro2}, or \cite{cof}) that the function $\varsigma$ satisfies
the conditions (\ref{supl}) (see Corollary \ref{csi} below). And hence, Theorem \ref{prin2} says that the Riemann Hypothesis is equivalent to $\rrp_\varsigma=\lp_0$.

\medskip

We guess that the converse of Theorem \ref{prin2} is also true:

\medskip

\noindent
\textit{Conjecture 1}. Let $B$ be a formal power series with $b_0=1$. If $\rrp_B=\lp_0$ then $B\in \lp I$, it is not a polynomial and
the limit (\ref{supl}) holds.

In fact, we have proved the following weak version of Conjecture 1. Write
$$
\rho_n=\frac{b_{n-2}b_n}{b_{n-1}^2}.
$$
\begin{Theo}\label{wic}
Let $B$ be a formal power series with $b_0=1$. If $\rrp_B=\lp_0$ then $B\in \lp I$, it is not a polynomial,
$$
\limsup_n\frac{b_{n-2}b_n}{b_{n-1}^2}=1,
$$
and if $(n_k)_k$ is an increasing sequence of positive integers such that $\lim_k\rho_{n_k}=1$, then $\lim_k\rho_{n_k-j}=1$ for all $j\ge 0$.
\end{Theo}

We have also proved Conjecture 1, under the additional assumption $\rrp_B=\rrp_B^p$.

In Section \ref{poz}, we prove that, under mild assumptions, if $A\in \rrp_B$, the real zeros of the Brenke polynomials $(p_n)_n$ generated by $A$ are simple and the zeros of $p_{n-1}$ interlace the zeros of $p_n$.

In Section \ref{ejps} we prove that the (generalized) hypergeometric functions
\begin{equation}\label{ihg}
B_\phi (z)=\pFq{0}{q}{-}{\phi_1,\dots,\phi_q}{z}, \quad \phi_i>0,
\end{equation}
belong to the \psa class and prove simplicity and interlacing properties of the zeros of the Brenke polynomials generated by $A\in \rrp_B$.

In the last section of this paper, we prove some more asymptotics for Brenke polynomials which provide new equivalencies for the Riemann Hypothesis.
Denote by $(q_{n,s})_n$ the Jensen polynomials generated by $\varsigma^{(s)}(z)/\gamma_s$ (see (\ref{var})):
\begin{equation}\label{nito}
q_{n,s}(z)=\frac{1}{\gamma_s}\sum_{j=0}^n\frac{\gamma_{s+j}}{(n-j)!\,j!}z^j,
\end{equation}
where $\gamma_n$ is given by (\ref{varg}). As we have already mentioned, RH is equivalent to the polynomials $q_{n,s}$, $n,s\ge 0$, having only real zeros.
Regarding this equivalency, the following Theorem  has been proved by Griffin, Ono, Rolen and Zagier (see \cite[Theorem 1]{GORZ}):

\begin{TheoA}\label{gorz} If $n\ge 1$, then $q_{n,s}$ has only real zeros for all sufficiently large $s$.
\end{TheoA}

\noindent
(See also \cite{gort,sull}).

Theorem \ref{gorz} is a consequence of a particular case of the following asymptotic for Brenke polynomials that we prove in Section \ref{ult}.

\begin{Cor}\label{ggorz} Let $A$ and $B$ be normalized formal power series satisfying:
\begin{equation}\label{cabc2}
b_n\not =0, \mbox{for all $n\ge 0$}, \quad \lim_n\frac{b_{n-1}b_{n+1}}{b_n^2}=1.
\end{equation}
Denote by $(q_{n,s})_s$ the Brenke polynomials generated by $A$ and associated to $\displaystyle \frac{B^{(s)}}{s!\,b_s}$. Then
\begin{equation}\label{asau}
\lim_s\frac{b_s}{(s+1)_nb_{n+s}}\left(\frac{(n+s+1)b_{n+s+1}}{b_{n+s}}z\right)^nq_{n,s}\left(\frac{b_{n+s}}{(n+s+1)b_{n+s+1}}\frac{1}{z}\right)=q_n(z),
\end{equation}
where $(q_n)_n$ are the Jensen polynomials generated by $A$.
Moreover, if $A\in \lp$ for $n\ge 0$ there exists $s_n$ such that the polynomial $p_{n,s}$ has only real zeros for $s\ge s_n$.
\end{Cor}

When $B=\varsigma$ (\ref{var}) and $A(z)=e^z$, we have that the polynomials $q_{n,s}$, $n,s\ge 0$, are the Jensen polynomials generated by $\varsigma^{(s)}(z)/\gamma_s$ (\ref{nito})
and then the asymptotic (\ref{asau}) gives
$$
\lim_s\frac{\gamma_s}{\gamma_{n+s}}\left(\frac{\gamma_{n+s+1}}{\gamma_{n+s}}z\right)^nq_{n,s}\left(\frac{\gamma_{n+s}}{\gamma_{n+s+1}}\frac{1}{z}\right)=\frac{1}{n!}(1+z)^n.
$$
From where Theorem \ref{gorz} by Griffin, Ono, Rolen and Zagier follows (compare with \cite[Corollary 3.2]{far}, see also \cite{Ki}).

Other asymptotics (that we prove in Section \ref{ult}) will lead to new equivalencies with the Riemann hypothesis.

\begin{Cor}\label{esf1} Let $\alpha$ be a real number $\alpha>-1$. The following are equivalent.
\begin{enumerate}
\item The Riemann Hypothesis is true.
\item The polynomial
$$
p_{n,s}^\alpha(z)=\frac{(-1)^s}{(\alpha+1)_s\gamma_s}\sum_{j=0}^n\frac{(\alpha+n-j+1)_s}{j!\,(n-j)!} \gamma_{n-j+s}z^j
$$
has only real zeros for all $n,s\ge 0$.
\end{enumerate}
Moreover, we have the asymptotic
$$
\lim_s\frac{(-1)^n\gamma_s}{\gamma_{n+s}(\alpha+s+1)_n}p_{n,s}^\alpha\left(-(\alpha+n+s+1)\frac{\gamma_{n+s+1}}{\gamma_{n+s}}z\right)=\frac{z^n}{(\alpha+1)_n}L_n^{\alpha}(1/z),
$$
where $L_n^\alpha$ is the $n$-th Laguerre polynomial (\cite[pp, 241-244]{KLS}).
Hence, for every $n\ge 0$ there exists $s_n$ such that the polynomial $p_{n,s}^\alpha$ has only real zeros for $s\ge s_n$.
\end{Cor}

Compare with the following \textit{dual} version.

\begin{Cor}\label{esf2}
Let $\alpha$ be a real number $\alpha>-1$. The following are equivalent.
\begin{enumerate}
\item The Riemann Hypothesis is true.
\item The polynomial
$$
q_{n,s}^\alpha(z)=\frac{1}{\gamma_s}\sum_{j=0}^n\frac{(-1)^{n-j}}{j!\,(n-j)!\,(\alpha+1)_{n-j}}\gamma_{j+s}z^j
$$
has only real zeros for all $n,s\ge 0$.
\end{enumerate}
Moreover, we have the asymptotic
$$
\lim_s\frac{\gamma_s}{\gamma_{n+s}}\left(\frac{\gamma_{n+s+1}}{\gamma_{n+s}}z\right)^nq_{n,s}^\alpha \left(\frac{\gamma_{n+s}}{\gamma_{n+s+1}}\frac{1}{z}\right)=\frac{1}{(\alpha+1)_n}L_n^\alpha (z).
$$
Hence, for every $n\ge 0$ there exists $s_n$ such that the polynomial $q_{n,s}^\alpha$ has only real zeros for $s\ge s_n$.
\end{Cor}

\section{Brenke polynomials}\label{bp}

Let $A$ and $B$ be two formal power series
$$
A(z)=\sum_{n=0}^\infty a_nz^n,\quad B(z)=\sum_{n=0}^\infty b_nz^n,
$$
normalized by taking $a_0=b_0=1$ (when that happens we refer to them as normalized formal power series).

Let us recall that the polynomials $(p_n)_n$ are Brenke polynomials generated by $A$ and associated to $B$ if
\begin{equation}\label{dgfx}
A(z)B(xz)=\sum_{n=0}^\infty p_n(x)z^n.
\end{equation}
(see \cite{Brenke}, \cite[Ch. V]{Ch} or \cite[p. 654]{Ism}).
The normalization $a_0=b_0=1$ implies that $p_0=1$ and that $p_n$ is a polynomial of degree less than or equal to $n$; more precisely, it is of degree $n$ if and only if $b_n\not=0$ and then the leading coefficient of $p_n$ is $b_n$.

The following expression for the Brenke polynomials generated by $A$ and associated to $B$ follows straightforwardly from (\ref{dgfx}):
\begin{equation}\label{bps}
p_{n}(x)=\sum_{j=0}^{n}a_jb_{n-j}x^{n-j}=\sum_{j=0}^{n}a_{n-j}b_{j}x^{j}.
\end{equation}

When $b_n\not =0$, $n\ge 0$, one can characterize Brenke polynomials exclusively from the formal power series $B$. Indeed,
we associate to $B$ the linear operator $\Lambda_B$ defined in the linear space of polynomials by
\begin{equation}\label{opgb}
\Lambda_B(x^n)=\begin{cases}0,& n=0,\\\displaystyle \frac{b_{n-1}}{b_{n}}x^{n-1}, &n\ge 1.\end{cases}
\end{equation}
It follows then easily that $(p_n)_n$ are Brenke polynomials associated to $B$ if and only if
\begin{equation}\label{opgb2}
\Lambda_B (p_n)(x)=p_{n-1}(x).
\end{equation}

\medskip

Many interesting sequences of polynomials are Brenke polynomials.
\begin{enumerate}
\item Appell polynomials: $B(z)=e^z$. Bernoulli, Euler or Hermite polynomials are examples of Appell polynomials.
\item Appell-Dunkl polynomials: $B(z)=E_\mu(z)$, where for $\mu\in \CC\setminus\{-1,-2,\dots \}$, we consider the entire functions
\begin{align}\nonumber
  \I_\mu(z) &= 2^\mu \Gamma(\mu+1) \frac{J_\mu(iz)}{(iz)^\mu}, \\\label{emu}
  \Ea(z) &= \I_\mu(z) + \frac{z}{2(\mu+1)} \, \I_{\mu+1}(z),
\end{align}
where $J_\mu$ is the Bessel function of order~$\mu$ (let us remark that $E_{-1/2}(z)=e^z$).
\item $q$-Appell. In two versions $B(z)=1/(z;q)_\infty$ (which it is often denoted by $e_q(z)$ and it is one of the two $q$-exponential functions), and $B(z)=(-z;q)_\infty$ (which it is often denoted by $E_q(z)$ and it is the other $q$-exponential function).
\end{enumerate}

\medskip

It is worth noticing that there is a symmetry between $A$ and $B$ in the definition of Brenke polynomials.

\begin{Lem}\label{sdc}
Let $(p_n)_n$ be the sequence of Brenke polynomials generated by $A$ and associated to $B$. Then, the polynomials $x^np_n(1/x)$ are the Brenke polynomials generated by $B$ and associated to $A$.
\end{Lem}

\begin{proof}
The proof is straightforward, because if
$$
A(z)B(xz)=\sum_{n=0}^\infty p_n(x)z^n,
$$
by setting $u=xz$ and $y=1/x$ then
$$
B(u)A(yu)=\sum_{n=0}^\infty y^np_n(1/y)u^n.
$$
\end{proof}

In the Appell case, $B(z)=e^z$, the symmetry in Lemma \ref{sdc} gives the Jensen polynomials (\ref{jepo}) generated by $A$.

\bigskip

Since we want to study Brenke polynomials with real zeros, along this paper we will assume that the formal power series $A$ and $B$ have real Taylor coefficients, that is, $a_n,b_n\in \RR$, for all $n\ge 0$.

\bigskip

The following Lemma will be also useful.

\begin{Lem}\label{ay} Let $B$ be a formal power series such that $b_3\not =0$ and $A\in \rrp_B$. If for some $l$, $a_l=a_{l+1}=0$  then $a_j=0$ for all $j\ge l$. \end{Lem}

\begin{proof}

The proof is a consequence of the following fact. If $p$ is a polynomial with only real zeros and $p'$ has at $\zeta\in \RR$ a zero of multiplicity bigger than 1, then $p(\zeta )=0$.

Indeed, since $A\in \rrp_B$, each polynomial $p_n$, $n\ge 0$, has only real zeros. Take $l$ the smallest positive integer such that $a_l=a_{l+1}=0$ (and hence $a_{l-1}\not =0$).
Using (\ref{bps}), we see that $p_{l+2}(x)=a_{l+2}+a_{l-1}b_3x^3+\cdots +b_{l+2}x^{l+2}$. Hence  $p_{l+2}'$ has degree at least $2$ and has at $0$ a zero of multiplicity bigger than $1$. Hence $a_{l+2}=p_{l+2}(0)=0$.
\end{proof}

\subsection{Asymptotic for Brenke polynomials}\label{abp}

In what follows, we assume that
\begin{equation}\label{bnn0}
b_n\not =0, \quad \mbox{for $n$ big enough}.
\end{equation}

The asymptotic will be proved under the assumption
\begin{equation}\label{pbsxx}
\lim_{n\to \infty}\frac{b_{n-2}b_n}{b_{n-1}^2}=1.
\end{equation}

We straightforwardly have the following Lemma.

\begin{Lem}\label{lems} Let $(b_n)_n$ be a sequence of numbers satisfying (\ref{bnn0}). Then the following conditions are equivalent.
\begin{enumerate}
\item The limit (\ref{pbsxx}) holds for the sequence $(b_n)_n$.
\item \begin{equation}\label{pbs}
\begin{cases} \mbox{there exists a sequence $(\tau_n)_n$, $\tau_n\not=0$, \mbox{$n$ big enough},}&  \\
\mbox{such that $\lim_{n\to \infty}\frac{b_{n-j}}{b_n\tau_n^j}=1$, $j\ge 1$}.&
\end{cases}
\end{equation}
\item There exists a sequence $(\tau_n)_n$, $\tau_n\not=0$, \mbox{$n$ big enough}, such that
\begin{equation}\label{pbsx}
\lim_{n\to \infty}\frac{b_{n-1}}{b_n\tau_n}=1,\quad \lim_{n\to \infty}\frac{\tau_{n-1}}{\tau_n}=1.
\end{equation}
\end{enumerate}
\end{Lem}
Let us note that a sequence $\tau_n$ satisfying both parts (2) and (3) in the previous Lemma is $\tau_n=b_n/b_{n+1}$.

\medskip

We are now ready to prove the asymptotic for Brenke polynomials stated  in Theorem \ref{iasy}.

\begin{proof}[Proof of Theorem \ref{iasy}]

Using (\ref{bps}) we write
\begin{equation}\label{las}
\left(\frac{z}{\tau_n}\right)^n\frac{p_{n}(\tau_n/z)}{b_n}=\sum_{j=0}^{n}a_j\frac{b_{n-j}}{b_n\tau_n^j}z^{j}.
\end{equation}

Hence, fixed $N$, we have for $n\ge N$,
$$
\left(\frac{z}{\tau_n}\right)^n\frac{p_{n}(\tau_n/z)}{b_n}-A(z)=\sum_{j=0}^{N}a_j\left(\frac{b_{n-j}}{b_n\tau_n^j}-1\right)z^{j}+
\sum_{j=N+1}^{n}a_j\left(\frac{b_{n-j}}{b_n\tau_n^j}-1\right)z^{j}+\sum_{j=n+1}^{\infty}a_jz^{j}.
$$
Since
$$
\left|\frac{b_{n-j}}{b_n\tau_n^j}-1\right|\le M+1,
$$
we get
$$
\left|\left(\frac{z}{\tau_n}\right)^n\frac{p_{n}(\tau_n/z)}{b_n}-A(z)\right|\le \left|\sum_{j=0}^{N}a_j\left(\frac{b_{n-j}}{b_n\tau_n^j}-1\right)z^{j}\right|+
(M+1)\sum_{j=N+1}^{\infty}|a_jz^{j}|.
$$
From where the theorem follows easily.
\end{proof}

\bigskip

Let us point out that if $B\in\lp I$ is not a polynomial, and we write
\begin{equation}\label{taue}
\tau_n=b_n/b_{n+1},
\end{equation}
then $B$ satisfies the assumption (\ref{cbt}) in Theorem \ref{iasy}.
Indeed, on the one hand, if we write
\begin{equation}\label{rho}
\rho_n=\frac{b_{n-2}b_n}{b_{n-1}^2},
\end{equation}
a simple computation shows that
\begin{equation}\label{eqnp}
\frac{b_{n-j}}{b_n\tau_n^j}=\prod_{i=0}^{j-1}\rho_{n+1-i}^{j-i}.
\end{equation}
On the other hand, it is well-known (see (\ref{blcon}) below) that if $\sum_{n=0}^\infty b_nz^n\in \lp I$ is not a polynomial, then the sequence $(b_n)_n$ is strictly log-concave, that is
$$
0<b_{n-2}b_n<b_{n-1}^2,\quad n\ge 2.
$$
Hence, we can conclude that $\rho_n<1$, and then (\ref{eqnp}) shows that $B$ satisfies the assumption (\ref{cbt}) in Theorem \ref{iasy}.

\bigskip

Write
$$
\varsigma (z)=\sum_{n=0}^\infty b_nz^n,
$$
where $\varsigma$ is the function defined by (\ref{var}).

\begin{Cor}\label{csi} The function $\varsigma$ (\ref{var}) satisfies that $b_n>0$ for all $n$ and
\begin{equation}\label{lisi}
\lim_{n\to \infty}\frac{b_{n-2}b_n}{b_{n-1}^2}=1.
\end{equation}
\end{Cor}

\begin{proof}
On the one hand,  let us consider the Riemann Xi function  $\Xi$, so that $\varsigma(z)=\Xi(-i\sqrt z)/\xi(1/2)$. Since $\Xi$ is an even function, we have
$$
b_n=\frac{(-1)^n\Xi ^{(2n)}(0)}{(2n)!\,\xi(1/2)}.
$$
The integral representation
$$
\Xi (x)=2\int_0^\infty \Phi(u)\cos (ux)du,
$$
where
$$
\Phi(u)=2\sum_{n=1}^\infty (2n^4\pi^2e^{9u/2}-3n^2\pi e^{5u/2})e^{-n^2\pi e^{2u}}\ge 0,\quad u>0
$$
(see \cite[(10.1.3) and (10.1.4)]{tit}) gives
$$
\Xi^{(2n)}(0)=2(-1)^n\int_0^\infty \Phi (u)u^ndu,
$$
from where easily follows that $b_n>0$.

On the other hand the limit (\ref{lisi}) is an easy consequence of the asymptotic for
$\Xi^{(2n)}(0)$ found by Grosswald (see (\cite[(12), p. 17]{gro1}, \cite{gro2} or \cite{cof}). If
$D_n=nb_n^2-(n+1)b_{n-1}b_{n+1}$ then $D_n=b_n^2(1+O(1/\log (n)))$.
Hence
$$
\left|\frac{D_n}{b_n^2}-1\right| =\left|n-(n+1)\frac{b_{n-1}b_{n+1}}{b_n^2}-1\right|\le \frac{M}{\log n},
$$
from where the limit (\ref{lisi}) follows easily.
\end{proof}

\bigskip

Theorem \ref{iasy} has the following dual version.

\begin{Theo}\label{iasya}
Assume that the sequence $(a_n)_n$, $a_0=1$, satisfies
\begin{equation}\label{pans}
\lim_{n\to \infty}\frac{a_{n-2}a_n}{a_{n-1}^2}=1,
\end{equation}
and let $(\mu_n)_n$ be a sequence satisfying
$$
\lim_{n\to \infty}\frac{a_{n-j}}{a_n\mu_n^j}=1,\quad \mbox{for all $j\ge 1$}.
$$
For an analytic function $B$ at $z=0$ (with radius of convergence $r$) consider the Brenke polynomials $(p_n)_n$ generated by $A(z)=\sum_na_nz^n$ and associated to $B$.
If there exist $M,N >0$ such that
\begin{equation}\label{cbta}
\left|\frac{a_{n-j}}{a_n\mu_n^j}\right|\le M,\quad j,n\ge N,
\end{equation}
then
\begin{equation}\label{las2a}
\lim_{n\to \infty}\frac{1}{a_n}p_{n}(z/\mu_n)=B(z).
\end{equation}
uniformly in compact sets of $\{z:\vert z\vert <r\}$.
\end{Theo}

\section{Setting the real-rooted Brenke polynomial problem and some illustrative examples}\label{rrbp}

Let us start with a definition.

\begin{Def}\label{rrpyp} Given a formal power series $B$, we say that the formal power series $A$ has the real-rooted polynomial property if  for all $n\ge 0$ the Brenke polynomial $p_n$ generated by $A$ has all its zeros real.
If there exists $n_0\ge 0$ such that  for $n\ge n_0$ the Brenke polynomial $p_n$ generated by $A$ has all its zeros real, we say that the formal power series $A$ has the almost real-rooted polynomial property. Finally, if the set $\{n: \mbox{the $n$-th Brenke polynomial $p_n$ has all its zeros real}\}$ is infinite, we say that the formal power series $A$ has the partial real-rooted polynomial property.
\end{Def}

For a formal power series $B$ we write
\begin{align}\label{rrpp}
\rrp _B&=\{\mbox{$A$: $A$ is a normalized formal power series }\\\nonumber
&\hspace{1.2cm}\mbox{having the real-rooted polynomial property}\},\\
\label{rrpp2}
\rrp _B^a&=\{\mbox{$A$: $A$ is a normalized formal power series }\\\nonumber
&\hspace{1.2cm}\mbox{having the almost real-rooted polynomial property}\},\\
\label{rrpp3}
\rrp _B^p&=\{\mbox{$A$: $A$ is a normalized formal power series }\\\nonumber
&\hspace{1.2cm}\mbox{having the partial real-rooted polynomial property}\}
\end{align}
(see (\ref{rrppi}) and (\ref{rrpp2i}) in the Introduction).

Obviously $\rrp_B\subset \rrp_B^a \subset \rrp_B^p$.

Several problems arise from the previous definition. The first one is the \textit{real-rooted Brenke polynomial problem}:
\begin{equation}\label{pob}
\mbox{Given a formal power series $B$, characterize the set $\rrp_B$.}
\end{equation}
And similarly, we have the \textit{almost and partial real-rooted Brenke polynomial problem} which consist in characterizing the sets $\rrp_B^a$ and $\rrp_B^p$, respectively.

We have not found these problems explicitly posed in the literature, although
they have been studied and solved for some particular $B$'s
and they are related to some other problems studied
during the last 100 years. Here it is some examples.

\subsection{The Appell case}
As explained in the Introduction, the solution for the case $B(z)=e^z$, i.e., the Appell case, goes more than one century back and it is due, independently, to Jensen \cite{jen}, and Pólya and Schur \cite{PS}. They proved that
$\rrp_{e^z}=\lp_0$ (see Theorem \ref{TPS}), where $\lp_0$ is the (normalized) Laguerre-Pólya class (see Definition \ref{def1} and (\ref{lpn})).

As stated in Theorem \ref{thiwi}, the functions in the Laguerre-Pólya class are those entire functions which can be approximate (locally uniformly) by polynomials with only real zeros.
As we wrote in the Introduction, Theorem \ref{thiwi} has the stronger version provided by Theorem \ref{tlli}, which we prove next.

\begin{proof}[Proof of Theorem \ref{tlli}]
Consider the Appell polynomials generated by the formal power series $A(z)=\sum_{j=0}^\infty a_{j}z^j$:
$$
p_n(z)=\sum_{j=0}^n\frac{a_{n-j}}{j!}z^j.
$$
The proof will follow if we prove that $p_n$ has only real zeros for all $n\ge 0$, because of Theorem \ref{TPS}.

Write $s_k$ for the polynomial defined by
$$
s_k(z)=z^{n_k}r_{n_k}(1/z)=\sum_{j=0}^{n_k}a_{n_k-j,n_k}z^j.
$$
The assumption implies that $s_k$ has only real zeros. And so, for $n\le n_k$, the polynomial $s_k^{(n_k-n)}(z)$ has only real zeros as well. Using that
$$
s_k^{(n_k-n)}(z)=\frac{n_k!}{n!}a_{0,n_k}z^n+\frac{(n_k-1)!}{(n-1)!}a_{1,n_k}z^{n-1}+\cdots +(n_k-n)!\,a_{n,n_k},
$$
we can conclude that the polynomial
$$
\frac{s_k^{(n_k-n)}(z/n_k)}{(n_k-n)!}=\frac{a_{0,n_k}}{n!}\frac{(n_k-n+1)_n}{n_k^n}z^n+\frac{a_{1,n_k}}{(n-1)!}\frac{(n_k-n+1)_{n-1}}{n_k^{n-1}}z^{n-1}+\cdots +a_{n,n_k}
$$
has only real zeros.

By taking limit when $k\to \infty$, we finally obtain that the polynomial
$$
\frac{a_{0}}{n!}z^n+\frac{a_{1}}{(n-1)!}z^{n-1}+\cdots +a_{n}=p_n(z)
$$
has only real zeros.
\end{proof}

\bigskip

The functions in the class $\lp I$ can be characterized as follows.

\begin{TheoA}\label{smac} An entire function $f(z)=\sum_{n=0}^\infty \theta_n z^n$ is of first type in the Laguerre-Pólya class if and only if is in the Laguerre-Pólya class and the sequence $(\theta_n)_n$ has constant or alternating sign.
\end{TheoA}

See \cite[p. 8]{CrCs}.

\medskip

The following definition will also be useful.

\begin{Def}\label{prr} Let $T$ be a linear operator acting in the linear space of polynomials. We say that $T$ preserves real-rootedness if for all polynomial $p$ having only real zeros then the polynomial $T(p)$
has only real zeros as well.
\end{Def}

The functions in the class $\lp I$ can also be characterized in terms of operators which preserve real-rootedness.

\begin{TheoA}\label{PS2} Given a sequence $\theta=(\theta_n)_n$ the linear operator
\begin{equation}\label{lops}
T_\theta (\sum_{j=0}^r\alpha_j x^j))=\sum_{j=0}^r\theta_j\alpha_jx^j
\end{equation}
preserves real-rootedness if and only if the power series
$$
A(z)=\sum_{n=0}^\infty \theta_n\frac{z^n}{n!}
$$
is of first type in the Laguerre-Pólya class.
\end{TheoA}

See \cite{PS}, \cite{CrCs}.

As a consequence we have.
\begin{TheoA}\label{llcon} Let $B$ be a function of first type in the Laguerre-Pólya class. Then
either $b_n\not=0$, $n\ge 0$, or there exits $n_0$, such that $b_n\not=0$, for $n\le n_0$, and $b_n=0$, for $n\ge n_0+1$.
Moreover, when $b_n\not=0$ then
\begin{equation}\label{blcon}
0\le b_{n-1}b_{n+1}<b_n^2,\quad n\ge 1
\end{equation}
(i.e. $(b_n)_n$ is strictly log-concave).
\end{TheoA}

See \cite[Lemma 3.10]{Bra}.

\medskip

\begin{Cor}\label{bessel} For a non negative integer $l$, the linear operator $T_{\theta^{[l]}}$ associated to the sequence $\theta^{[l]}=(1/(n+l)!)_n$
preserves real-rotedness.
\end{Cor}

\begin{proof}
A simple computation shows
$$
\sum_{n=0}^\infty \frac{z^n}{(n+l)!\,n!}=\frac{J_l(2\sqrt{-z})}{(-z)^{l/2}},
$$
where $J_l$ is the Bessel function of order $l$.

Since $\frac{J_l(2\sqrt{-z})}{(-z)^{l/2}}$ is an entire function of order $1/2$ and positive zeros, we deduce that $\frac{J_l(2\sqrt{-z})}{(-z)^{l/2}}\in \lp I$. Hence, the corollary is an easy consequence of Theorem~\ref{PS2}.
\end{proof}

\bigskip

For the Appell case, we have that $\rrp_{e^z}=\rrp_{e^z}^a=\rrp_{e^z}^p$. Indeed, Appell polynomials are characterize because $p_n'=p_{n-1}$, and since the derivative operator $d/dx$ preserves real-rootedness, we trivially have that $\rrp_{e^z}=\rrp_{e^z}^a=\rrp_{e^z}^p$.

We point out that for a formal power series $B$, the operator $\Lambda_{B}$ (see (\ref{opgb})) plays the role of $d/dx$ for $e^z$. This suggests the following definition.

\begin{Def}\label{sdco} We say that a formal power series $B$  is stable if the operator $\Lambda_B$ (\ref{opgb}) preserves real-rootedness (see Definition \ref{prr}).
\end{Def}

As a consequence of (\ref{opgb2}) we have the following corollary.

\begin{Cor}\label{feo} Assume that the formal power series $B$ is stable and that for certain $n$, $p_n$ has only real zeros. Then for $0\le j\le n-1$, the  polynomial $p_j$ has also only real zeros.
In particular $\rrp_{B}=\rrp_B^a=\rrp_{B}^p$.
\end{Cor}

Write $\theta=(\theta_n)_{n=0}^\infty$ for the sequence
$$
\theta_n=\begin{cases} 0, &n=0,\\ \frac{b_{n-1}}{b_n},&n\ge 1.\end{cases}
$$
Since $\Lambda_B=\frac{1}{x}T_\theta$, Theorem \ref{PS2} gives:

\begin{TheoA}\label{PS2x} The formal power series $B$ is stable if and only if the formal power series
\begin{equation}\label{elp}
C(z)=\sum_{n=0}^\infty \frac{b_{n}}{b_{n+1}}\frac{z^n}{(n+1)!}
\end{equation}
 is an entire function of first type in the Laguerre-Pólya class.
\end{TheoA}

We can extend the linear operator $\Lambda _B$  from the linear space of polynomials to that of formal power series
as follows:
$$
\Lambda_B\left(\sum_{n=0}^\infty a_nz^n\right)=\sum_{n=0}^\infty a_{n+1}\frac{b_n}{b_{n+1}}z^n.
$$
We then have the following corollary.

\begin{Cor}\label{sccd} If $B$ is stable, then the Laguerre-Pólya class is closed under the operator $\Lambda_B$. Moreover, $\lp I$ is also close under the operator $\Lambda_B$.
\end{Cor}

\begin{proof}

Indeed, if $A\in \lp$ then there exists a sequence of polynomials $(r_n)_n$ with only real zeros such that $\lim_nr_n(z)=A(z)$, uniformly in $\{z: |z|\le r\}$ for some $r>0$ (Theorem \ref{thiwi}).
Since $\Lambda_B$ is stable, we have that the polynomials
$\Lambda_B(r_n)$ have only real zeros. The proof follows now by applying Theorem \ref{tlli}.

Using Theorems \ref{PS2x} and \ref{smac}, we deduce that if $B$ is stable, then the sequence $(b_{n}/b_{n+1})_n$ has constant or alternating sign. Hence, proceeding as before we can also prove that if $A\in \lp I$ then
$\Lambda_B(A)\in \lp I$.
\end{proof}

\medskip

In the case $B=\varsigma$ (\ref{var}), we have checked that the operator $\Lambda_\varsigma$ does not preserve real-rootedness. Indeed, the polynomial $\Lambda_\varsigma ((x+1)^4)$ has degree $3$ but only one real zero.

\subsection{Entire functions with real-rooted Taylor sections}\label{tayr}
Consider next the following example $B(z)=1/(1-z)$. Since
$$
\sum_{j=0}^na_{n-j}x^j=x^{n}\sum_{j=0}^na_{j}/x^j,
$$
we straightforwardly deduce that $A\in \rrp_B$ if and only if the polynomial
$\sum_{j=0}^na_{j}x^j$ has only real zeros for all $n\ge 0$.

Hence the problem of characterizing $\rrp_B$ for $B(z)=1/(1-z)$ is equivalent
to the problem of characterizing the analytic functions whose Taylor sections have only real zeros. This is a well-known problem with a long tradition of more than one century (see, for instance, \cite{Ost}, and the references therein).

We next consider an example related to the previous one which shows that also when $B$ is a formal power series converging only for $z=0$ the characterization of $\rrp_B$ could be an interesting problem. Indeed,  for $q>1$, let $B_q$ be the formal power series
\begin{equation}\label{sfbq}
B_q(z)=\sum_{n=0}^\infty q^{n^2}z^n.
\end{equation}
For $q=1$, we have $B_1(z)=1/(1-z)$.

For a formal power series $A$, a simple computation gives
$$
p_n(x)=\sum_{j=0}^na_jq^{(n-j)^2}x^{n-j}=
q^{n^2}x^n\sum_{j=0}^na_jq^{j^2}\frac{1}{(q^{2n}x)^{j}}.
$$
Hence $A\in \rrp_{B_q}$ if and only if the polynomial $\sum_{j=0}^na_jq^{j^2}x^{j}$ has only real zeros for all $n\ge 0$. Hence
\begin{equation}\label{ciq}
\rrp_{B_q}=\left\{\sum_{n=0}^\infty a_nz^n: \sum_{n=0}^\infty a_n q^{n^2}z^n\in \rrp_{B_1}\right\}.
\end{equation}
Since
\begin{equation}\label{cagi}
\sum_{n=0}^\infty \frac{z^n}{n!\,q^{n^2}}\in \lp I
\end{equation}
(see \cite{KLV2}), we have that $\rrp_{B_q}\subset \lp$. Indeed, if $A\in \rrp_{B_q}$, then (\ref{ciq}) says that $\sum_{j=0}^na_jq^{j^2}x^{j}$ has only real zeros. Hence Theorem \ref{PS2} implies that
$\sum_{j=0}^na_jx^{j}$ has only real zeros (because (\ref{cagi})), and then using Theorem \ref{thiwi} we can conclude that $A\subset \lp$.

For $s>1$, take now the functions
$$
A_{s,1}(z)=\sum_{n=0}^\infty \frac{1}{s^{n^2}}z^n,\quad A_{s,2}(z)=\sum_{n=0}^\infty \frac{1}{n!\,s^{n^2}}z^n.
$$
The problem of whether $A_{s,i}\in \rrp_{B_q}$ has already been considered in the literature (although using other terminology). Indeed, according to (\ref{ciq}), $A_{s,1}\in \rrp_{B_q}$ if and only if $\sum_{j=0}^n\frac{1}{(s/q)^{j^2}}z^j$ has only real zeros for all $n\ge 0$. Similarly, $A_{s,2}\in \rrp_{B_q}$ if and only if $\sum_{j=0}^n\frac{1}{j!\,(s/q)^{j^2}}z^j$ has only real zeros for all $n\ge 0$. These problems were solved
by Katkova, Lobova and Vishnyakova in \cite{KLV} and \cite{KLV2}, respectively:
\begin{enumerate}
\item $\sum_{j=0}^n\frac{1}{(s/q)^{j^2}}z^j$ has only real zeros for all $n\ge 0$ if and only if $s\ge 2q$.
\item There exists a constant $q_\infty\approx 3,23\dots$, such that $\sum_{j=0}^n\frac{1}{j!\,(s/q)^{j^2}}z^j$ has only real zeros for all $n\ge 0$ if and only if $s^2\ge q^2q_\infty$.
\end{enumerate}
In particular, $\{A_{s,1},s\ge 2q\}\cup \{A_{s,2}, s^2\ge q^2q_\infty\}\subset \rrp_{B_q}$.

\subsection{An example where $\rrp_B$ is trivial.}

We next show how to use the asymptotic properties of the sequence $(b_n)_n$ to characterize $\rrp_B$. This illustrates that to use asymptotic properties of the sequence $(b_n)_n$ (as stated
in Theorems \ref{iasy} and \ref{iasya}) will be a useful tool to study the real-rooted Brenke polynomial problem of characterizing the set $\rrp_B$ (\ref{rrppi}).

Consider the rational function
$$
B(z)=\frac{1}{1-z^3}+\frac{-z+z^2}{1-z^3/2}=\sum_{j=0}^\infty b_jz^j,
$$
where
\begin{equation}\label{epm}
b_n=\begin{cases}1,& \mbox{if $n=3k$},\\
-2^{-k},& \mbox{if $n=3k+1$},\\
2^{-k},& \mbox{if $n=3k+2$}.\end{cases}
\end{equation}
We will next show that the set $\rrp_B$ is trivial, i.e. $\rrp_B=\{ az+1: a\in \RR\}.$

Indeed, if $A(z)=\sum_{j=0}^\infty a_{j}z^{j}\in \rrp_B$, we have that the polynomials
$$
z^{3n}p_{3n}(1/z)=\sum_{j=0}^{3n}a_jb_{3n-j}z^j
$$
have only real zeros. Since
$$
\lim_na_jb_{3n-j}=\begin{cases} a_j, &j=3s,\\0, &\mbox{otherwise},\end{cases}
$$
using Theorem \ref{tlli}, we conclude that
$$
C(z)=\sum_{j=0}^\infty a_{3j}z^{3j}\in \lp.
$$
In particular, this shows that the function $\tilde C(z)=C(z^{1/3})$ is entire of order at most $2/3$. Assume next that $\tilde C$ is not the constant function. Hence it has complex zeros, and then the function $C$ has to have non real zeros. But this contradicts that $C\in \lp$. Hence $\tilde C$ has to be the constant function, and so we deduce that $a_{3j}=0$, $j\ge 1$.

Proceeding in a similar way using the polynomials
$$
z^{3n+i}p_{3n+i}(1/z),\quad i=1,2,
$$
we deduce that also $a_{3j+i}=0$, $j\ge 1$. And then
$A(z)=1+a_1z+a_2z^2$. But since
$$
p_n(x)=b_nx^{n-1}(x^2+a_1b_{n-1}x/b_n+a_2b_{n-2}/b_n),
$$
it is easy to conclude from (\ref{epm}) that $a_2=0$.

\section{The \psa class}\label{jps}
The examples in the last section show that the solution of the problem of characterizing the set $\rrp_B$  may strongly depend on the formal power series $B$.

In what follows we will concentrate in the following problem. We say that a formal power series $B$ belongs to the \psa class if $\rrp_B=\lp$.

The initials \textbf{\psa} stand for \textbf{J}ensen, \textbf{P}olya and \textbf{S}chur, because we want to characterize the formal power series $B$ such that $\rrp_B=\rrp_{e^z}=\lp$,  and the characterization of $\rrp_{e^z}$ is due to Jensen, Pólya and Schur.

In the next subsections, we will find necessary and sufficient conditions on a formal power series to be in the \psa class.

\subsection{A necessary condition for a formal power series $B$ to be in the \psa class}

We start proving a stronger version of Theorem \ref{bpoli}, which characterizes the class $\lp I$ in terms of Brenke polynomials with real zeros. In particular, this Theorem implies that if a formal power series $B$ belongs to the \psa class then $B\in \lp I$.

\begin{Theo} \label{bpol}
Let $B$ be a formal power series with $b_0=1$. The following conditions are equivalent.
\begin{enumerate}
\item $\{e^z,1-z^2,(1+z)^l,l\ge 2\} \subset \rrp_B$.
\item $B$ is an entire function of the first type in the Laguerre-Pólya class.
\item For all $l\ge 0$, $\sum_{j=0}^\infty \frac{b_j}{(j+1)_lj!}z^j\in \lp I$.
\item $\lp_0 \subset \rrp _B$.
\end{enumerate}
\end{Theo}

Before proving the theorem we need the following lemma.

\begin{Lem}\label{lgi}
Let $B$ be a normalized formal power series  satisfying that
$$
\{1-z^2,(1+z)^l,l\ge 2\} \subset \rrp_B.
$$
Then one of the following conditions holds
\begin{enumerate}
\item $b_n\not =0$ for all $n$.
\item $B$ is a polinomial of degree $k$ and $b_n\not =0$, $n=0,\dots , k$.
\end{enumerate}
Write $\kappa=\begin{cases} \infty, &\mbox{if $B$ satisfies (1)},\\ k+1,&\mbox{if $B$ satisfies (2)}.\end{cases}$. Then for $2\le n< \kappa $, $b_{n}b_{n-2}>0$, and
\begin{equation}\label{coti}
\frac{b_{n-1}^2}{b_{n-2}b_{n}}\ge 1+\frac{1}{n^2-1}> 1.
\end{equation}
\end{Lem}

\begin{proof}[Proof of Lemma \ref{lgi}]

For $A(z)=1-z^2\in\rrp_B$, a simple computation using (\ref{bps}) gives
$$
p_n(x)=x^{n-2}\left(b_nx^2-b_{n-2}\right).
$$
Since $p_n$ has to have only real zeros, we deduce that $b_{n-2}b_n\ge 0$ for all $n\ge 2$.

If $b_n\not =0$ for all $n\ge 0$ then $b_{n-2}b_n> 0$ for all $n\ge 2$.

If there exists $n\ge 2$ such that $b_n=0$, let $k$ be the smallest positive integer for which $b_{k+1}=0$. We next prove that $b_j=0$ for $j\ge k+1$. Indeed, since $(z+1)^2\in \rrp_B$, a simple computation using (\ref{bps}) gives that the polynomial
$$
p_n(x)=x^{n-2}\left(b_nx^2+2b_{n-1}x+b_{n-2}\right)
$$
has to have only real zeros for all $n\ge 2$. Setting $n=k+2$ we have
$$
p_{k+2}(x)=x^{k}\left(b_{k+2}x^2+b_{k}\right),
$$
and since $b_{k+2}b_{k}\ge 0$, $b_k\not=0$, we deduce that $b_{k+2}=0$.

Assume then that $b_j=0$, for $j=k+1,\dots , k+l-1$ and $l\ge 3$. Since $(z+1)^{l}\in \rrp_B$, a simple computation using (\ref{bps}) gives that
$$
p_n(x)=x^{n-l}\sum_{s=0}^{l} \binom{l}{s}b_{n-s}x^{l-s}
$$
has to have only real zeros for all $n\ge 2$. Setting $n=k+l$ we have
$$
p_{k+l}(x)=x^{k}\left(b_{k+l}x^{l}+b_{k}\right),
$$
and we deduce that $b_{k+l}=0$.

This proves that $B(z)=\sum_{j=0}^{k}b_jx^j$, with $b_j\not=0$, $j=0,\dots, k$, and $b_{n-2}b_n> 0$ for all $2\le n<k+1$.

Consider finally a formal power series $A(z)=\sum_{n=0}^\infty a_nz^n$, so that (see (\ref{bps}))
$$
p_n(x)=\sum_{j=0}^na_{n-j}b_jx^j.
$$
Assume that $p_n$ has  only real zeros. Taking $(n-2)$-th derivatives, we deduce that also the polynomial
$$
\frac{n!}{2}b_nz^2+(n-1)!\,a_1b_{n-1}z+(n-2)!\,a_2b_{n-2}
$$
has also only real zeros. This gives that
\begin{equation}\label{cho}
\frac{b_{n-1}^2}{b_{n-2}b_n}\ge 2\left(1+\frac{1}{n-1}\right) \frac{a_2}{a_1^2}.
\end{equation}
The inequality (\ref{coti}) follows by taking $A(z)=(1+z)^{n+1}$.
\end{proof}

Let us remark that the inequality (\ref{coti}) implies that
\begin{equation}\label{bee}
\begin{cases} \mbox{$b_{n-1}/b_n$, $1\le n< \kappa $, is an increasing sequence},& \mbox{if $b_1>0$},\\
\mbox{$b_{n-1}/b_n$, $1\le n< \kappa $, is a decreasing sequence},& \mbox{if $b_1<0$}.
\end{cases}
\end{equation}

\begin{proof}[Proof of Theorem \ref{bpol}]

We firstly prove (1) $\Rightarrow$ (2). Let $A(z)=e^z$, and write $(p_n)_n$ for the Brenke polynomials generated by $A$ associated to $B$. Write $\tilde p_n(x)=x^np_n(1/x)$. Since $A\in \rrp_B$, the polynomial $p_n$ has only real zeros for all $n\ge 0$, and then $\tilde p_n$ has also only real zeros. Lemma \ref{sdc} says that $\tilde p_n$ are the Brenke polynomials generated by $B$ associated to $e^z$, that is,  $\tilde p_n$ are the Appell polynomials generated by $B$. Hence, $B\in \rrp_{e^z}$. According to  Theorem \ref{TPS} $\rrp_{e^z}=\lp$ and then $B$ has to be an entire function in the Laguerre-Pólya class.

Since
$$
\{1-z^2,(1+z)^l,l\ge 2\} \subset \rrp_B,
$$
Lemma \ref{lgi} shows that either $b_{n-2}b_n>0$, $2\le n< +\infty$ or $b_{n-2}b_n>0$, $2\le n< k+1$ (if $B$ is a polynomial of degree $k$). In any case this shows that
the sequence $(b_n)_n$ has equal or alternating sign and then $B\in \lp I$.

\bigskip

We next prove (2) $\Rightarrow$ (3).
Consider the linear operator $T_B$ acting in the linear space of polynomials and defined by $T_B(x^j)=j!\,b_jx^j$.
Using Theorem \ref{PS2}, we have that $T_B$ preserves real-rootedness. Corollary \ref{bessel} says that also the
linear operator $T_{\theta^{[l]}}$, $l$ a nonnegative integer, acting in the linear space of polynomials and defined by $T_{\theta^{[l]}}(x^j)=x^j/(j+l)!$
preserve real-rootedness.
And so, the operator $T_{B,l}$ defined by
\begin{equation}\label{tbl}
T_{B,l}(x^j)=T_{\theta^{[l]}}\circ T_B(x^j)=\frac{j!\,b_jx^j}{(j+l)!}=
\frac{b_jx^j}{(j+1)_l}
\end{equation}
preserves real-rootedness as well. Using again Theorem \ref{PS2}, we deduce that
$$
\sum_{j=0}^\infty \frac{b_j}{(j+1)_lj!}x^j\in \lp I.
$$

\bigskip

The proof of (3) $\Rightarrow$ (4) is as follows. Using Theorem \ref{PS2}, the assumption (3) is equivalent to assume that the operators $T_{B,l}$ (\ref{tbl}), $l\ge 0$, preserve real-rootedness.

Let $A\in\lp_0$.  Assume first that $A$ is a polynomial with real zeros. Write $k$ for the degree of $A$. On the one hand we have for $n\ge k$
$$
p_n(x)=\sum_{j=0}^ka_jb_{n-j}x^{n-j}
=\sum_{j=n-k}^na_{n-j}b_{j}x^{j}=T_{B,0}\left( \sum_{j=n-k}^na_{n-j}x^{j}\right).
$$
And on the other hand
$$
\sum_{j=n-k}^na_{n-j}x^{j}=x^nA(1/x).
$$
Obviously the polynomial $x^nA(1/x)$ has also only real zeros and then $p_n$ does, because $T_{B,0}$ preserves real-rootedness.

For $n=k-l$, with $0\le l\le k-2$, we first derive $l$ times the polynomial $\tilde A(x)=x^kA(1/x)=\sum_{j=0}^ka_{k-j}x^j$:
$$
\tilde A^{(l)}(x)=\sum_{j=l}^{k}j(j-1)\cdots (j-l+1)a_{k-j}x^{j-l}=\sum_{j=0}^{k-l}(j+1)_la_{k-l-j}x^j.
$$
Since $\tilde A$ has only real zeros, the derivative $\tilde A^{(l)}$ has only real zeros as well. The definition (\ref{tbl}) shows
$$
T_{B,l}(\tilde A^{(l)})(x)=\sum_{j=0}^{k-l}a_{k-l-j}b_jx^j=p_{k-l}(x).
$$
Hence $p_{k-l}$ has real zeros because $T_{B,l}$ preserves real-rootedness.

Assume finally that $A$ is not a polynomial. Since $A$ is in the Laguerre-Pólya class, there exists a sequence of polynomials
$(A_N)_N$,  all their zeros are real and $\lim_NA_N(z)=A(z)$ uniformly in $\{z:|z|\le r\}$ for some $r>0$ (Theorem \ref{thiwi}). Denote by $(p_{n,N})_n$ the sequence of Brenke polynomials generated by $A_N$.
We have already proved that $p_{n,N}$ has only real zeros.

It is clear that fixed $n$, we have
$$
p_n(z)=\lim_{N\to \infty}p_{n,N}(z)
$$
uniformly in compact sets of $\CC$. Since the zeros of each $p_{n,N}$ are all real, we conclude that all the zeros of $p_n$ has to be  real as well.

\bigskip

The proof of (4) $\Rightarrow$ (1) is straight forward because $\{e^z,1-z^2,(1+z)^l,l\ge 2\} \subset \lp_0$.
\end{proof}

We next prove Corollary \ref{rhe}, which provides an equivalency of the Riemann Hypothesis in terms of real-rootedness of a sequence of polynomials.

\begin{proof}[Proof of Corollary \ref{rhe}]

RH is equivalent to $\varsigma\in \lp I$ (\ref{var}). Theorem \ref{bpol} says that $\varsigma\in \lp I$ if and only if for all $l\ge 0$,
$$
\sum_{j=0}^\infty \frac{\gamma_j}{(j+1)_lj!^2}z^j\in \lp I.
$$
Repeating the process we have that for each positive integer $N$, RH is equivalent to
$$
\varsigma_{l_1,\dots,l_N}(z)=\sum_{j=0}^\infty \frac{\gamma_j}{j!\,\prod_{i=1}^N(j+l_i)!}z^j\in \lp I
$$
for all non-negative integers $l_i$, $i=1,\dots, N$.
Since $\gamma_j>0$, $j\ge 0$ (Corollary \ref{csi}), we deduce that for all non-negative integers $l_i$, $i=1,\dots, N$, $\varsigma_{l_1,\dots,l_N}(z)\in \lp I$ if and only if $\varsigma_{l_1,\dots,l_N}(z)\in \lp $.
And this is equivalent to the Jensen polynomials generated by $\varsigma_{l_1,\dots,l_N}$ and defined as
$$
q_{l_1,\dots,l_N;n}(x)=\sum_{j=0}^n \frac{\gamma_j}{j!\,(n-j)!\,\prod_{i=1}^N(j+l_i)!}x^j
$$
having only real zeros for $n\ge 0$ (because of Theorem \ref{TPS}). We have that $\hat q_{N,n}=q_{n,\dots,n;n}$ (see (\ref{jenpx})).

Hence RH implies that $\hat q_{N,n}$ has only real zeros for all $n\ge 0$.

Assume next that $\hat q_{N,n}$ has only real zeros for all $n\ge 0$.
A simple computation shows that
$$
(x^n\hat q_{N,n}(x))^{(n-l_1)}=x^{l_1}q_{l_1,n,\dots,n;n}(x)
$$
for $n\ge l_1$. And so $q_{l_1,n,\dots,n;n}$ has only real zeros for all $n\ge l_1$. Repeating the process, we conclude that $q_{l_1,\dots,l_N;n}$ has only real zeros for all $n$ and
$l_i\le n$, $i=1,\dots, N$.

Write $r_{l_1,\dots,l_N;n}(x)=x^nq_{l_1,\dots,l_N;n}(1/x)$, so that $r_{l_1,\dots,l_N;n}$ has only real zeros for all $l_i,n\ge 0$ if and only if $q_{l_1,\dots,l_N;n}$ has only real zeros for all $l_i, n\ge 0$. Hence, we have already proved that $r_{l_1,\dots,l_N,n}$ has real zeros for all $n$ and $l_i\le n$. If $\{i: l_i\ge n\}=\{i_0\}$,
a simple computation shows that
$$
r_{l_1,\dots,l_N;l_{i_0}}^{(l_{i_0}-n)}(x)=r_{l_1,\dots,l_N;n}(x)
$$
for $0\le n\le l_{i_0}$. From where we deduce that also $q_{l_1,\dots,l_N;n}$ has only real zeros for all $n$ when $l_i\le n$ except for one $i$.
The proof can be completed by iterating the process.
\end{proof}

\medskip

We conclude this section with a couple of more results.

The first one is a characterization of the polynomials in $\rrp_B^a$ under the assumption (\ref{pbsxx}). It is a consequence of Theorem \ref{iasy}. Indeed, if $A$ is a polynomial with $k=\deg A$, then (\ref{las}) gives
$$
\frac{1}{\tau_n^n}\frac{p_{n}(\tau_nx)}{b_n}=x^{n-k}\sum_{j=0}^{k}\frac{b_{n-j}}{b_n\tau_n^{j}}a_{j}x^{k-j}.
$$
Hence, if $\tau_n$ satisfies (\ref{pbs}), we deduce that
\begin{equation}\label{lasf}
\lim_{n\to \infty}\frac{1}{\tau_n^n}\frac{p_{n}(\tau_nx)}{b_nx^{n-k}}=x^{k}A(1/x)
\end{equation}
uniformly in compact sets of $\CC$.

\begin{Cor}\label{pfc}
Let $A$ and $B$ be a polynomial and a normalized formal power series satisfying (\ref{pbsxx}), respectively.
\begin{enumerate}
\item If $A\in \rrp_B^a$, then all the zeros of $A$ has to be real.
\item If  all the zeros of $A$ are real and simple then $A\in \rrp_B^a$.
\item Assume in addition that $B\in \lp_0$ and that $\{n:b_{n-2}b_n\le 0\}$ is a finite set. If all the zeros of $A$ are real then $A\in \rrp_B^a$.
\end{enumerate}
\end{Cor}

\begin{proof}
The first part is a straightforward consequence of the uniform convergence in compact sets of (\ref{lasf}).

We next prove the second part. Since $A(0)=a_0=1$, we have $A(0)\not =0$. Hence the polynomial  $\tilde A(z)=z^{k}A(1/z)$ is also a polynomial of degree $k$ whose zeros are real and simple.

Since $p_n$ has at $x=0$ a zero of multiplicity  $n-k$, it is enough to prove that $p_n$ has $k$ real zeros different to $0$.
Denote
$$
\zeta_{1}<\cdots <\zeta_k
$$
for the $k$ zeros of $\tilde A$. Write $D_i=\{z:\vert z-\zeta_i\vert <\epsilon\}$, $i\in \{1,\dots, k\}$, and choose $\epsilon>0$ such that $D_i\cap D_j=\emptyset$, $i\not =j$, and $0\not \in D_i$. Using  Hurwitz Theorem, we deduce that for $n$ big enough each polynomial $p_{n}(\tau_nx)$ has exactly one zero in each disc $D_i$. Since $p_{n}(\tau_nx)$ has real coefficients, we conclude that those zeros has to be real and different to $0$. Hence, $p_{n}$ has $k$ real zeros different to $0$, and a zero of multiplicity at least $n-k$ at $x=0$. Hence, all the zeros of $p_n$ have to be real.

The proof of the third part is as follows. Write $m-1$ for the maximum element of $\{n:b_{n-2}b_n\le 0\}$ and take $m$-derivative of the function $B$
$$
B^{(m)}(z)=\sum_{n=0}^\infty (n+1)_mb_{n+m}z^n.
$$
Since $B^{(m)}\in \lp$ and $b_{n+m}b_{n+m+2}>0$, $n\ge 0$, we conclude that $B^{(m)}\in \lp I$. And hence $T_{B^{(m)},m}$ preserves real-rootedness (where $T_{B^{(m)},m}$ is the operator (\ref{tbl}) for $B^{(m)}$).
If the polynomial $A$ has degree $k$ and only real zeros, proceeding as in the proof of (3) $\Rightarrow$ (4) in Theorem \ref{bpol} we have for $n\ge m+k$
\begin{align*}
p_n(x)&=\sum_{j=0}^na_jb_{n-j}x^{n-j}=\sum_{j=n-k}^na_{n-j}b_{j}x^{j}=\sum_{j=n-k}^na_{n-j}\frac{(j-m+1)_m}{(j-m+1)_m}b_{j}x^{j}\\ &=x^mT_{B^{(m)},m}\left(\sum_{j=n-m-k}^{n-m}a_{n-m-j}x^{j}\right).
\end{align*}
We then deduce that $p_n$ has only real zeros because
$$
\sum_{j=n-m-k}^{n-m}a_{n-m-j}x^{j}=x^{n-m}A(1/x)
$$
has only real zeros and $T_{B^{(m)},m}$ preserves real-rootedness.
\end{proof}

In the second part of Corollary \ref{pfc}, the hypothesis of $A$ having simple zeros can not be removed. Here it is a counterexample: $B(z)=1+\sum_{n=1}^\infty z^n/n$, and $A(z)=(x-1)^2$.
A simple computation gives, for $n\ge 3$,
$$
p_n(x)=\frac{x^{n-2}}{n}\left(x^2-\frac{2n}{n-1}x+\frac{n}{n-2}\right).
$$
And it is easy to check that $p_n$ has always two zeros that are not real.

Theorem \ref{bpol} and Corollary \ref{pfc} imply that if a formal power series $B$ satisfies (\ref{pbsxx}) and it is stable, then $B\in \lp I$.

\begin{Theo}\label{rec}
Let $B$ be a normalized formal power series satisfying (\ref{pbsxx}).
Assume in addition that $B$ is stable. Then $\lp_0\subset \rrp_B$. As a consequence, $B\in \lp I$.
\end{Theo}

\begin{proof}
Assume first that $A$ is a polynomial with simple zeros. Since $A$ is in the Laguerre-Pólya class, $A$ only has real zeros. The second part of Corollary \ref{pfc} and Corollary \ref{feo} imply that $A\in \rrp_B^a=\rrp_B$.

Assume next that $A$ is a polynomial with multiple zeros. For a given $\epsilon>0$ we can modify each multiple zero to construct a polynomial $A_\epsilon$ whose zeros are all real and simple and such that $\lim_{\epsilon \to 0}A_\epsilon (z)=A(z)$. Denote $p_{n,\epsilon}$ for the Brenke polynomials generated by $A_\epsilon$. On the one hand, $A_\epsilon$ is obviously also in the Laguerre-Pólya class and hence for all $n$ $p_{n,\epsilon}$ has only real zeros.
On the other hand, it is straightforward to see that for $n\ge 0$, also $\lim_{\epsilon \to 0}p_{n,\epsilon} (z)=p_n(z)$. Since all the zeros of $p_{n,\epsilon}$ are real, so are the zeros of $p_n$.

Assume finally that $A$ is not a polynomial. Since $A$ is in the Laguerre-Pólya class, there exists a sequence of polynomials
$(A_N)_N$,  all their zeros are real and $\lim_NA_N=A$.
Hence we have already proved that $A_N\in \rrp_B$.
Denote by $(p_{n,N})_n$ the sequence of Brenke polynomials generated by $A_N$.
It is clear that fixed $n$, we have
$$
p_n(z)=\lim_{N\to \infty}p_{n,N}(z)
$$
uniformly in compact sets of $\CC$. Since the zeros of each $p_{n,N}$ are all real, we conclude that all the zeros of $p_n$ has to be  real as well.
\end{proof}

\medskip

\subsection{Sufficient conditions for a formal power series $B$ to be in the \psa class}\label{cssp}

We start proving that if $B$ is a polynomial, then $B$ is not in the \psa class.

\begin{Lem}\label{cfp}
Let $B$ be a polynomial of degree $k$ satisfying $b_0=1$, then $\rrp_B\not =\lp_0$. More precisely,
\begin{enumerate}
\item If for some $n_0$, $2\le n_0\le k$, $b_{n_0-1}^2-b_{n_0-2}b_{n_0}\le 0$ then $\lp_0\not\subset \rrp_B$.
\item If for all $n$, $2\le n\le k$, $b_{n-1}^2-b_{n-2}b_n>0$, then $\rrp_B\not\subset\lp_0$.
\end{enumerate}
\end{Lem}

\begin{proof}
The first part of the lemma is a consequence of the second part of Lemma~\ref{lgi}.

The second part of the lemma will follow if we prove that there is a polynomial $A\in \rrp_B$ of degree 2, $a_0=1$, without real zeros, and hence $A\not \in \lp_0$.

If we write $A(z)=1+a_1z+a_2z^2$, then, for $2\le n$
$$
p_n(x)=x^{n-2}(a_2b_{n-2}+a_1b_{n-1}x+b_nx^2).
$$
In order to prove that $A\in \rrp_B$ it is enough to prove that $p_n$ has only real zeros for $2\le n\le k$. But this is equivalent to the inequalities
\begin{equation}\label{inch}
a^2_1b_{n-1}^2-4a_2b_{n-2}b_n\ge 0,\quad 2\le n\le k.
\end{equation}
Take $a_2=a_1^2/4+\epsilon$, where $\epsilon>0$. This implies that $A$ has no real zeros. Substituting in (\ref{inch}), we get
$$
a^2_1b_{n-1}^2-4a_2b_{n-2}b_n=a^2_1(b_{n-1}^2-b_{n-2}b_n)-4\epsilon b_{n-2}b_n.
$$
Since $b_{n-1}^2-b_{n-2}b_n>0$, for $2\le n\le k$, by taking $\epsilon>0$ small enough
we see that the inequalities (\ref{inch}) hold.
\end{proof}

\medskip

We are now ready to prove Theorem \ref{prin2}.

\begin{proof}[Proof of Theorem \ref{prin2}]

On the one hand, since $B\in \lp I$, using Theorem \ref{bpol} we deduce that $\lp_0\subset \rrp_B$.

On the other hand, take a formal power series $A\in \rrp_B^p$ and write $(p_n)_n$ for the Brenke polynomials generated by $A$. Then there exits an increasing sequence $(n_k)$ of positive integers such that
the polynomials $(p_{n_k})_k$ have only real zeros. Since the limit (\ref{pbsxx}) holds, we can take a sequence $(\tau_n)_n$ satisfying (\ref{pbs}).
And, so we deduce that the polynomials
$$
r_{n_k}(z)=\left(\frac{z}{\tau_{n_k}}\right)^{n_k}\frac{p_{n_k}(\tau_{n_k}/z)}{b_{n_k}}=\sum_{j=0}^{n_k}a_j\frac{b_{n_k-j}}{b_{n_k}\tau_{n_k}^j}z^{j},
$$
have only real zeros (because $p_{n_k}$ have only real zeros). Using the limits (\ref{pbs}) we have
$$
\lim_ka_j\frac{b_{n_k-j}}{b_{n_k}\tau_{n_k}^j}=a_j.
$$
By applying Theorem \ref{tlli}, we deduce that $A\in \lp$, and so $\rrp_B^p\subset \lp_0$.

Hence
$$
\lp_0\subset \rrp_B\subset \rrp_B^a\subset \rrp_B^p\subset \lp_0.
$$
That is $\rrp_B= \rrp_B^a= \rrp_B^p= \lp_0$.
\end{proof}

Let us remark that we have actually proved the following.

\begin{Cor}\label{prinito} Assume $B\in \lp I$, $B(0)=1$, it is not a polynomial and
$$
\lim_n \frac{b_{n-2}b_n}{b_{n-1}^2}=1.
$$
Then $\rrp_B=\rrp_B^a=\rrp_B^p=\lp_0$.
\end{Cor}

\bigskip

In order to prove Theorem \ref{wic}, we need the two following lemmas.

\begin{Lem}\label{cflm}
Let $B$ be a formal power series with $b_n\not =0$, $n\ge n_0$, for some $n_0$,  and $b_0=1$. Assume also that there exists $0<\lambda<1$ and an increasing sequence $(n_k)_k$ of positive integers, $n_0\ge 2$, such that
\begin{equation}\label{lm1}
\frac{ b_{n_k-2}b_{n_k}}{b_{n_k-1}^2}\le\lambda,\quad k\ge 0.
\end{equation}
Then $\rrp_B^p\not\subset \lp_0$.
In particular, if for certain $s\ge 0$, $n_k=k+s$, we have $\rrp_B^a\not \subset \lp_0$, and $\rrp_B\not \subset \lp_0$ if $s=0$.
\end{Lem}

\begin{proof}
We prove that there is a polynomial $A\in \rrp_B^p$ of degree 2, $a_0=1$, without real zeros, and hence $A\not \in \lp_0$.

If we write $A(z)=1+a_1z+a_2z^2$, and proceed as in Lemma \ref{cfp}, we conclude that
$A\in \rrp_B^p$ follows if we prove that
\begin{equation}\label{inch2}
a^2_1b_{n_k-1}^2-4a_2b_{n_k-2}b_{n_k}\ge 0,\quad k\ge 2.
\end{equation}
Take $a_2=a_1^2/4+\epsilon$, where $\epsilon>0$. This implies that $A$ has no real zeros. Substituting in (\ref{inch2}), we get
$$
a^2_1b_{n_k-1}^2-4a_2b_{n_k-2}b_{n_k}=b_{n_k-1}^2\left[a^2_1\left(1-\frac{b_{n_k-2}b_{n_k}}{ b_{n_k-1}^2}\right)-\epsilon \frac{b_{n_k-2}b_{n_k}}{ b_{n_k-1}^2}\right].
$$
Using (\ref{lm1}), we deduce that for all $k\ge 0$
$$
a^2_1\left(1-\frac{b_{n_k-2}b_{n_k}}{ b_{n_k-1}^2}\right)-\epsilon \frac{b_{n_k-2}b_{n_k}}{ b_{n_k-1}^2}\ge a_1^2(1-\lambda)-\epsilon.
$$
Hence, by taking $\epsilon>0$ small enough we see that the inequalities (\ref{inch2}) hold.
\end{proof}

Consider again the notation (\ref{rho}) and write
\begin{equation}\label{ro2}
\rho_n=\frac{ b_{n-2}b_{n}}{b_{n-1}^2}.
\end{equation}

\begin{Lem}\label{ule}
Let $B$ be an analytic function at $z=0$ with $b_n\not =0$, $n\ge n_0$ and $b_0=1$.
Assume that
$$
0<\rho_n\le 1,\quad n\ge n_0+1,
\quad\quad
\limsup_n \rho_n=1,
$$
and $(z+1)^3\in \rrp_B^a$. If $(n_k)_k$ is an increasing sequence of positive integers
such that
$\lim_k \rho_{n_k}=1,$
then
$\lim_k \rho_{n_k-j}=1$, for all $j\ge 0$.
\end{Lem}

\begin{proof}
We proceed by reductio to absurdum, and assume that there exists an increasing sequence of positive integers $(n_k)_k$
such that
\begin{equation}\label{ldc}
\lim_k \rho_{n_k}=1,\quad \mbox{and \quad $\lim_k \rho_{n_k-1}=\lambda<1$.}
\end{equation}
Writing $\tau_n=b_{n}/b_{n+1}$, the identity (\ref{eqnp}) gives
\begin{equation}\label{inm}
\frac{b_{n-j}}{b_{n}\tau_{n}^j}=\rho_{n+2-j}\rho_{n+3-j}^2\cdots \rho_{n+1}^{j}.
\end{equation}
Let $(p_n)_n$ be the Brenke polynomials generated by $(z+1)^3$ associated to $B$. Since $(z+1)^3\in\rrp_B^a$, we deduce that
the polynomial
$$
\left(\frac{z}{\tau_n}\right)^n\frac{p_{n}(\tau_n/z)}{b_n}=z^{n-3}\left(1+3\frac{b_{n-1}}{b_n\tau_n}z+3\frac{b_{n-2}}{b_n\tau_n^2}z^{2}+\frac{b_{n-3}}{b_n\tau_n^3}z^{3}\right)
$$
has only real zeros for $n$ big enough. Using (\ref{inm}), this polynomial can be rewritten in the form
$$
z^{n-3}\left(1+3(\rho_{n+1}z)+3\rho_{n}(\rho_{n+1}z)^{2}+\rho_{n-1}\rho_{n}^2(\rho_{n+1}z)^{3}\right).
$$
Hence, we deduce that the polynomial
$$
1+3z+3\rho_{n}z^{2}+\rho_{n-1}\rho_{n}^2z^{3}
$$
has also real zeros for $n$ big enough.

By setting $n=n_k$, taking limit as $k\to \infty$ and using (\ref{ldc}) we conclude
that the polynomial
$$
1+3z+ 3z^2+\lambda z^3
$$
has only real zeros as well. And so the polynomial
$$
\lambda +3z+ 3z^2+z^3
$$
has also only real zeros. But this is a contradiction because that polynomial has the two following non-real zeros
$$
(1-\lambda)^{1/3}-1-\frac{\sqrt 3(1-\lambda)^{1/3}}{2}i,\quad (1-\lambda)^{1/3}-1+\frac{\sqrt 3(1-\lambda)^{1/3}}{2}i.
$$
\end{proof}

We are now ready to prove Theorem \ref{wic} which it is a weaker version of Conjecture~1.

\begin{proof}[Proof of Theorem \ref{wic}]

On the one hand, Theorem \ref{bpol} and Lemma \ref{cfp} implies that $B\in \lp I$ and it is not a polynomial, and then $0<\rho_n<1$ (see Theorem \ref{llcon}), where
$\rho_n$ is defined by (\ref{ro2}).

On the other hand, Lemma \ref{cflm} implies that $\limsup_n \rho_n=1$, and since $(1+z)^3\in \lp_0=\rrp_B\subset \rrp_B^a$, Lemma \ref{ule} shows that if $(n_k)_k$ is an increasing sequence of positive integers such that
$\lim_k \rho_{n_k}= 1$ then $\lim_k \rho_{n_k-j}= 1$ for all $j\ge 0$.
\end{proof}

Theorem \ref{wic} implies the following.

\begin{Lem} Let $B\in \lp_0 I$ which it is not a polynomial. Define
$$
\mathcal X=\{a\in[0,1]: \mbox{$a$ is a limit point of $(\rho_n)_n$}\}.
$$
If $\{1\}\subsetneq \mathcal X$,  and $1$ is an isolated point of $\mathcal X$ then $\lp_0\subsetneq\rrp_B$.
\end{Lem}

\begin{proof}
Take $\epsilon >0$ such that $\mathcal X\cap [1-\epsilon,1]=\{1\}$. Define
$$
X=\{n: 1-\epsilon <\rho_n<1\}.
$$
Since $1\in \mathcal X$, we have that $X$ is infinite and $\lim_{n\in X}\rho_n=1$ (because $1$ is an isolated point of $\mathcal X$). Moreover, if $Y\subset \NN$ is infinite and $X\cap Y$ is finite then $\lim_{n\in Y}\rho_n\not =1$ (in the sense that either the limit $\lim_{n\in Y}\rho_n$ does not exist or if there exists is different to $1$).
We next prove that $Y=(-1+X)\setminus X$ is infinite. Indeed, assume on the contrary that $(-1+X)\setminus X$ is finite.  Hence, there exists a positive integer $l$ such that
if $n\in -1+X$ and $n\ge l$ then $n\in X$, that is, if $m\in X$ and $m\ge l+1$ then $m-1\in X$. Take $n_0\in X$ such that $n_0\ge l+1$, and given a positive number $M>0$ take $k\ge M$ such that
$k+n_0\in X$ (such $k$ always exists because $X$ is infinite). Since $-i+k+n_0\ge x_0\ge l$, for $i=1,\dots k$, we deduce that
$-i+k+n_0\in X$. As a consequence, we have that $\{n_0,n_0+1,n_0+2,\dots\}\subset X$. But this would imply $\lim_n\rho_n=1$ which it contradicts  $\{1\}\subsetneq \mathcal X$.

Since $Y=(-1+X)\setminus X$ is infinite, we can take an increasing sequence of positive integers $n_k-1\in Y$ (and so $n_k\in X$), such that the limit  $\lim_k\rho_{n_k-1}$ exists. Since $Y\cap X=\emptyset$,
we deduce that $\lim_k\rho_{n_k-1}=\lambda<1$, and, since $n_k\in X $, $\lim_k\rho_{n_k}=1$. Theorem \ref{wic} then implies $\lp_0\subsetneq\rrp_B$.
\end{proof}

\medskip

We next prove Conjecture 1 under the additional hypothesis $\rrp_B=\rrp_B^p=\lp_0$.

\begin{Cor}
If $\rrp_B^p=\lp_0$ then $B\in \lp $, $\{n:b_{n-2}b_n\le 0\}$ is a finite set and the limit
(\ref{supl}) holds.
Moreover, if $\rrp_B=\rrp_B^p=\lp_0$, then $B\in \lp I$, it is not a polynomial and the limit
(\ref{supl}) holds.
\end{Cor}

\begin{proof}
If $\rrp_B^p=\lp_0$, then $e^z\in \rrp_B^p$. Using Lemma \ref{sdc}, we deduce that $B\in \rrp_{e^z}^p=\rrp_{e^z}=\lp_0$.

If there exists an increasing sequence $(n_k)_k$ of positive integer, such that $b_{n_k-2}b_{n_k}\le 0$ then $A(z)=1+z+5z^2\in \rrp_B^p$. Indeed, we have $p_n(z)=z^{n-2}(5b_{n-2}+b_{n-1}z+b_nz^2)$.
Since $b_{n_k-2}b_{n_k}\le 0$, we deduce that $p_{n_k}$ has only real zeros and then $A\in \rrp_B^p$, which contradicts $\rrp_B^p=\lp_0$ because $A\not \in \lp_0$.
And so  $\{n:b_{n-2}b_n\le 0\}$ has to be a finite set.

Proceeding as in the last part of the proof of Lemma \ref{lgi}, we deduce that if $A\in \rrp_B$ then for $n$ big enough
$$
\frac{b_{n-1}^2}{b_{n-2}b_n}\ge 2\left(1+\frac{1}{n-1}\right) \frac{a_2}{a_1^2}.
$$
Using that $B \in \rrp_{e^z}$, we deduce that $e^z \in \rrp_{B}$ (Lemma \ref{sdc}), and hence $\rho_n<1$ for $n$ big enough. And so using Lemma \ref{cflm} we deduce that $\lim_n\rho_n=1$.

Assume next that $\rrp_B=\rrp_B^p=\lp_0$. Theorem \ref{wic} implies then that $B\in \lp I$ and it is not a polynomial.
Using the first part of the corollary we deduce that $\lim_n\rho_n=1$.
\end{proof}

\section{Properties of zeros of Brenke polynomials}\label{poz}
For two normalized formal power series $A$ and $B$ with $A\in \rrp_B$, the Brenke polynomials $(p_n)_n$ generated by $A$ and associated to $B$ has only real zeros.
In this Section we study the properties of the zeros of $(p_n)_n$. More precisely, simplicity and interlacing properties.

\begin{Def}\label{aqv} Let $q$ and $p$ be two polynomials with only real zeros and $\deg q=k$, $\deg p=k+1$. Denote by $\zeta_1\le \dots \le \zeta_k$ the real zeros of $q$ and by $\eta_1\le \dots \le \eta_{k+1}$ the real zeros of $q$.
\begin{enumerate}
\item We say that the zeros of  $q$  interlace the zeros of $p$ if
\begin{equation}\label{intz}
\eta_1\le \zeta_1\le \eta_2 \le \dots \le \eta_{k-1}\le \zeta_k\le \eta_{k+1}.
\end{equation}
\item If all the above inequalities are strict, we say that the zeros of $q$ strictly interlace the zeros of $p$.
\item If $x=\lambda$ is a common zero of $p$ and $q$, with multiplicities $l+1$ and $l$, respectively, we say that the zeros of  $q$  strictly interlace the zeros of $p$ except for the common zero $x=\lambda$ if
all the inequalities in (\ref{intz}) are strict except for $\eta_i=\zeta_i=\dots =\zeta_{i+l-1}=\eta_{i+l}=\lambda$.
\end{enumerate}
\end{Def}

\medskip

In general, the zeros of the Brenke polynomials generated by $A\in\rrp_B$ do not have to be simple or interlaced, as the following example shows.
Consider
$$
B(z)=\sum_{n=0}^\infty 2^{n^2}z^n,\quad A(z)=\sum_{n=0}^\infty z^n/4^{n^2}.
$$
Then, using (\ref{bps}), we have
$$
p_n(z)=2^{n^2}z^n\sum_{j=0}^n\frac{1}{2^{j^2}}\frac{1}{(4z)^j}.
$$
Since each polynomial $\sum_{j=0}^nz^j/2^{j^2}$ has only real zeros (see \cite{Huc} or, also, \cite{KLV}), we conclude that $A\in\rrp_B$. But $p_2(z)=(64z+1)^2/256$ has a double zero which it is smaller than the zeros of $p_3(z)=(512z+1)(262144z^2+1536z+1)/262144$.

\bigskip

However, under mild conditions we can prove interlacing properties for the zeros of $(p_n)_n$.

\begin{Lem}\label{enz} Let $B$ a normalized formal power series satisfying that if $A\in \rrp_B$ then $(1+bz)A\in\rrp_B$, for all $b\in \RR$.
Then for $A\in \rrp_B$, the zeros of $p_{n-1}$ interlaces the zeros of $p_n$. In particular, if $\zeta$ is a zero of multiplicity $k$ of $p_n$ then it is a zero of multiplicity $k-1$ of $p_{n-1}$.
\end{Lem}

\begin{proof}

The lemma is a consequence of Obreshkov Theorem \cite[Theorem 8]{Branden}

\begin{TheoA}\label{obre}
Let $p$ and $q$ be polynomials with only real zeros and with $\deg (p)=\deg(q)+1$. Then the following are equivalent.
\begin{enumerate}
\item The zeros of $p$ (strictly) interlace the zeros of $q$.
\item For each real numbers $\alpha,\beta$ the polynomial $\alpha p+\beta q$ has only real (and simple) zeros.
\end{enumerate}
\end{TheoA}

Indeed, in order to prove that the zeros of $p_{n-1}$ interlaces the zeros of $p_n$, it will be enough to show that for $b\in \RR$, the polynomial $p_n+bp_{n-1}$ has only real zeros.
But $p_n+bp_{n-1}$, $n\ge 1$, are the Brenke polynomials generated by the function $(1+bz)A\in \rrp_B $, and so the zeros of $p_n+bp_{n-1}$ are real.
\end{proof}

The assumption in Lemma \ref{enz} are always satisfied by the functions in the \psa class.

\bigskip

Simplicity of the zeros is more demanding. For instance, it is not true even in the Appell case. Indeed, consider $B(z)=e^z$ and $A(z)=e^{az}$, $a\in \RR$. Then, the Appell polynomials defined by $A$ are $p_n(x)=(x+a)^n$, and they have a multiple zero of order $n$ at $x=a$. Surprisingly enough the (quasi) simplicity of the zeros holds if we assume that the only multiple zero is $x=0$.

\begin{Theo}\label{tal} Let $B$ be a normalized formal power series with $b_1,b_3\not =0$, $n\ge 0$. Consider the Brenke polynomials generated by the function $A\in \rrp_B$, and assume that if $\zeta $ is a zero of $p_n$ of multiplicity larger than $1$ then $\zeta=0$.
Then
\begin{enumerate}
\item If $A$ is a polynomial of degree $k$ then for $0\le n\le k$, $p_n$ has simple zeros. For $n\ge k+1$, we have
\begin{equation}\label{pesfxg}
p_n(x)=x^{n-k}r_{n}(x),
\end{equation}
where $r_n$ is a polynomial of degree $k$ with simple zeros.
\item If $A$ is not a polynomial then $p_n$ has simple zeros.
\end{enumerate}
\end{Theo}

\begin{proof}

The proof of the first part is as follows.
Using (\ref{bps}), we deduce (\ref{pesfxg}). For $n\ge k$, we also deduce that
$r_{n}(0)=a_k\not=0$. Since the only multiple zero has to be $0$, this implies that $r_n$
has simple zeros. For $0\le n<k$, we proceed by \textit{reductio ad absurdum}. Indeed, we have $p_n(0)=a_n$ and $p_n'(0)=a_{n-1}b_1$.
If $p_n$ has a multiple zero at $x=\zeta$, since $\zeta=0$, we deduce that $a_n=a_{n-1}=0$. Since $n<k$, using Lemma \ref{ay}, we deduce that $a_k=0$, which it is a contradiction.

If $A$ is not a polynomial, then $p_n(0)=a_n$ and $p_n'(0)=a_{n-1}b_1$, and we can proceed as before using again Lemma \ref{ay}.
\end{proof}

\section{Examples of entire functions in the \psa class}\label{ejps}

The  following hypergeometric functions belong to the \psa class.

\begin{Cor}\label{ejh} If $p\le q$, $c_j>0$, $j=1,\dots, q$, and $m_j\in \NN\setminus \{0\}$, $j=1,\dots, p$, the (generalized) hypergeometric functions
\begin{equation}\label{mhe}
B(z)=\pFq{p}{q}{c_1+m_1,\dots,c_p+m_p}{c_1,\dots,c_q}{z}
\end{equation}
belong to the \psa class, i.e.,  $\rrp_B=\lp_0$.
\end{Cor}

\begin{proof}
On the one hand, it was proved in \cite{KaKa} (see also \cite{Ric}) that $B\in \lp I$. Then Theorem \ref{bpol} shows that $\lp_0\subset \rrp_B$.

On the other hand, it is easy to check that the limit (\ref{supl}) holds, Hence, we can apply  Theorem \ref{prin2} to deduce that $\rrp_B\subset \lp_0$.
\end{proof}

\bigskip

The case $p=0$ is specially interesting and will be studied in detail.

\bigskip

\begin{Def} Given a nonnegative integer $N$, let $\phi=(\phi_i)_{i=1}^N$ where $\phi_i>0$, $i=1,\dots N$. Associated to $\phi$ we define the hypergeometric function
\begin{equation}\label{h0f}
B_\phi (z)=\pFq{0}{N}{-}{\phi_1,\dots,\phi_N}{z}=\sum_{j=0}^\infty \frac{z^j}{j!\,\prod_{i=1}^N(\phi_i)_j}.
\end{equation}
\end{Def}

Hurwitz proved that $B_\phi\in \lp I$:

\begin{TheoA}\label{thu} For $\phi_i>0$, the hypergeometric function $\pFq{0}{N}{-}{\phi_1,\dots,\phi_N}{z}$ is of first type in the Laguerre-Pólya class.
\end{TheoA}

See \cite{hur} (also \cite{sok}).

Corollary \ref{ejh}   gives that for $\phi_i>0$ then $\rrp _{B_{\phi}}=\lp_0$.

Let us notice that for $N=1$ the Brenke polynomials generated by $A(z)=e^{z}$ and associated to $B_\phi$ are the Laguerre polynomials $(L_n^{\phi-1}(-x)/(\phi)_n)_n$.

The following theorem shows that only when $p=0$ the hypergeometric functions $B$ (\ref{mhe})  are stable.

\begin{Theo}\label{hedc} If $p\le q$, $c_j>0$, $j=1,\dots, q$, and $m_j\in \NN\setminus \{0\}$, $j=1,\dots, p$, then the hypergeometric function $B$ (\ref{mhe}) is stable if and only if $p=0$.
\end{Theo}

\begin{proof}
Using Theorem \ref{PS2x}, we have to characterize when the entire function
$$
C(z)=\sum_{n=0}^\infty \frac{b_{n}}{(n+1)!\,b_{n+1}}z^{n}
$$
belongs to $\lp I$.

A simple computation gives
\begin{equation}\label{cqs1}
\frac{b_{n}}{(n+1)b_{n+1}}=\frac{\prod_{i=1}^q(n+c_i)}{\prod_{i=1}^p(n+c_i+m_i)}.
\end{equation}

Using
$$
n+c=\frac{c(1+c)_n}{(c)_n},
$$
we have that
$$
C(z)=\frac{\prod_{i=1}^q c_i}{\prod_{i=1}^p(c_i+m_i)}\pFq{p+q}{p+q}{c_1+1,\dots,c_q+1,c_1+m_1,\dots,c_p+m_p}{c_1,\dots,c_q,c_1+m_1+1,\dots,c_p+m_p+1}{z}.
$$
Using \cite[Theorem 3]{KiKim} we can conclude that $C\in \lp I$ if and only if $p=0$.
This complete the proof.
\end{proof}

\bigskip

We next prove that for $A\in \rrp_{B_\phi}$ the zeros of the Brenke polynomials $(p_n)_n$ generated by $A$ are simple (with the only exception of $x=0$ when $A$ is a polynomial) and the zeros of $p_{n-1}$ strictly interlace with the zeros of $p_n$ (again with the only exception of $x=0$ when $A$ is a polynomial).

\begin{Theo}\label{pzgf0} Let $B_\phi$ be the hypergeometric function (\ref{h0f}),
where $\phi_i$ are positive real numbers. Consider $A\in \rrp_{B_\phi}=\lp_0$ and the Brenke polynomials $(p_n)_n$ generated by $A$.
\begin{enumerate}
\item If $A$ is a polynomial of degree $k$ then for $0\le n\le k$, $p_n$ has simple zeros and the zeros of $p_{n-1}$ strictly interlace the zeros of $p_{n}$. For $n\ge k$, we have
\begin{equation}\label{pesfxj}
p_n(x)=x^{n-k}r_{n}(x)
\end{equation}
where $r_n$ is a polynomial of degree $k$ with $r_n(0)\not =0$ and simple zeros, and the zeros of $p_n$ strictly interlace the zeros of $p_{n+1}$ except for the common zero $x=0$ (see Definition \ref{aqv}).
\item If $A$ is not a polynomial then $p_n$ has simple zeros and the zeros of $p_n$ strictly interlace the zeros of $p_{n+1}$.
\end{enumerate}
\end{Theo}

We need the following lemma.

\medskip

\begin{Lem}\label{jod} Assume $\phi_i>0$ and consider $A\in \rrp_{B_\phi}$ and the Brenke polynomials $(p_n)_n$ generated by $A$.
If $\zeta$ is a multiple zero of $p_n$ then $\zeta =0$.
\end{Lem}

\begin{proof}

We start with the following remark.
\begin{Rem}\label{urdc}
Let $B$ be a normalized formal power series such that $b_n\not =0$, $n\ge 0$. Take a normalized formal power series $A$ and  the Brenke polynomials $p_n$ (\ref{bps}) generated by $A$. Given $\theta_0\not =0$ consider the linear operator
$\Upsilon_B$ acting on polynomials as follows:
$$
\Upsilon_B(x^j)=\begin{cases} \theta_0,&j=0,\\
\frac{b_{j-1}}{jb_j}x^j,&j\ge 1.\end{cases}
$$
It is easy to check that
$$
\Lambda_B=\frac{d}{dx}\Upsilon _B,
$$
where $\Lambda_{B}$ is the operator defined in (\ref{opgb}) associated to $B$.

If we define $q_n=\Upsilon_B p_n/\theta_0$. A simple computation shows that the polynomials $q_n$ are also Brenke polynomials generated by $A$ but now associated to
\begin{equation}\label{tib}
\tilde B(z)=1+\sum_{j\ge 1}\frac{b_{j-1}}{\theta_0 j}z^j.
\end{equation}
Notice that when it makes sense $\tilde B'=B/\theta_0$.
\end{Rem}

In the case of the hypergeometric functions $B_\phi$, it is easy to see that, for $\theta_0=\prod_{i=1}^N(\phi_i-1)$,
\begin{equation}\label{oplh2}
\Upsilon_{B_\phi}=\prod_{i=1}^ND_{\phi_i-1},
\end{equation}
where for $\alpha\in \RR$, $D_\alpha$ denotes the first order differential operator
\begin{equation}\label{opa}
D_\alpha=\alpha I+x\frac{d}{dx}.
\end{equation}
Proceeding as in the proof of Theorem \ref{hedc}, it is easy to prove that for $\alpha_i\ge 0$, $i=1,\dots, N$, the differential operator
\begin{equation}\label{opax}
D=\prod_{i=1}^ND_{\alpha_i}
\end{equation}
preserves real-rootedness.

We next proceed in three steps.

\medskip

\noindent
\textit{Step 1}. Let $\alpha$ and $p$ be  a non-negative real number and a polynomial of degree $k$ which has only real zeros, respectively.
Assume that the polynomial  $D_\alpha p$ ($D_\alpha$ defined by (\ref{opa})) has a zero at $x=\zeta\not =0$ of multiplicity $l>1$. Then $p$ has a zero at $x=\zeta$ of multiplicity $l+1$.

\begin{proof}[Proof of Step 1]
Write $q=D_\alpha p=\alpha p+xp'$. Since $\alpha\ge 0$, $q$ has also degree $k$. We then have that all the zeros of $p$ and $q$ are real (because $D_\alpha$ preserves real-rootedness) and $q$ has a multiple zero at $\zeta$.

Assume first that all the zeros of $p$ are simple. Hence, $p'$ has $k-1$ simple zeros $\zeta_1<\dots <\zeta_{k-1}$. Write $\zeta_0=-\infty$ and $\zeta_k=+\infty$. For a polynomial $\lambda$ we define $\lambda(\pm\infty)=\sign (\lim_{x\to \pm \infty}p(x))$, where $\sign (x)=\begin{cases}1,& x>0,\\ -1,& x<0.\end{cases}$.

On the one hand $p(\zeta_i)p(\zeta_{i+1})<0$, $i=0,\dots, k-1$, and hence $q(\zeta_i)q(\zeta_{i+1})<0$ as well. Since $q$ has degree $k$, that shows that  $q$ has exactly one cero in each interval $(\zeta_i,\zeta_{i+1})$, $i=0,\dots, k-1$. Hence, the zeros of $q$ are also simple. Since this contradicts the assumption, we deduce that $p$ has to have at least a multiple zero.

Let $u_i$, $i=1,\dots , w$, the zeros of $p$, $u_i\not =0$, of multiplicity $v_i>1$. All of them are zeros of $p'$ of multiplicity $v_i-1$. Since $q=\alpha p+xp'$, $u_i$ is also a zero of $q$ of multiplicity $v_i-1$. We also have that $0$ is a zero of $p$ of multiplicity  $h$ if and only if it is a zero of $q$ of the same multiplicity $h$.
Define next the polynomial $r$ as
$$
r(x)=x^{h}\prod_{i=1}^w(x-u_i)^{v_i-1}
$$
(where $h=0$ if $p(0)\not =0$).
Write $\hat p=p/r$, $\hat q=q/r$ and $s=p'/r$ (if $h>0$, then $s$ is a rational function with a simple pole at $x=0$). By construction $\hat p=p/r$ is a polynomial with simple zeros, degree $m$  and $\hat p=\hat q+xs$. Write $\zeta_1<\dots <\zeta_{m-1}$ for the zeros of $s$. They are simple. Moreover, they are the zeros of $p'$ that are not zeros of $p$. Define as before $\zeta_0=-\infty$ and $\zeta_m=+\infty$. Given $\zeta_i$, $i=0,\dots, m-1$, write $X_i=\{j:\zeta_i<u_j<\zeta_{i+1}\}$ and $n_{X_i}$ for the number of elements in the set $X_i$. Hence $n_{X_i}$ is the number of multiple zeros of $p$ in $(\zeta_i,\zeta_{i+1})$. Write finally $m_i$ for the number of simple zeros of $p$ in  $(\zeta_i,\zeta_{i+1})$. Since $\zeta_i$, $i=1,\dots, m-1$, are zeros of $p'$, we get that there is at least one zero of $p$ in $(\zeta_i,\zeta_{i+1})$, $i=0,\dots, m-1$, and so $n_{X_i}+m_i\ge 1$. Now, on the one hand, the number of zeros of $p$ in $(\zeta_i,\zeta_{i+1})$ counting multiplicities is equal to $m_i+\sum_{j\in X_i}v_j$; and, on the other hand, the number of zeros of $r$ in $(\zeta_i,\zeta_{i+1})$ is equal to $\sum_{j\in X_i}(v_j-1)$. Since $\hat p=p/r$, we deduce that the number of zeros of $\hat p$ in
$(\zeta_i,\zeta_{i+1})$, $i=0,\dots, m-1$, is equal to $n_{X_i}+m_i\ge 1$. We then deduce that $\hat p$ has exactly one zero in each interval $(\zeta_i,\zeta_{i+1})$, $i=0,\dots, m-1$ (because $\hat p$ has degree $m$).
So, $\hat p(\zeta_i)\hat p(\zeta_{i+1})<0$, $i=0,\dots, m-1$. This gives $\hat q(\zeta_i)\hat q(\zeta_{i+1})<0$, $i=0,\dots, m-1$ (because $\hat p(\zeta_i)=\hat q(\zeta_i)$).
We conclude that $\hat q$ has simple zeros, and so $\zeta$ has to be equal to $u_i$ for some $i=1,\dots , w-1$. This implies that the multiplicity of $\zeta$ as a zero of $p$ is $l+1$.
\end{proof}

\medskip

\noindent
\textit{Step 2}. Consider the differential operator $D$ defined in (\ref{opax}), where $\alpha_i\ge 0$, $i=1,\dots, N$. If $p$ has only real zeros and $\zeta\not =0$ is a zero of $D p$ of multiplicity $l>1$, then $p$ has at $\zeta$ a zero of multiplicity $l+N$.

\begin{proof}[Proof of Step 2]

Obviously, Step 1 is the case $N=1$. The general case can be proved easily as follows.
If we write $\hat \alpha=(\alpha_2,\dots, \alpha_N)$, then the differential operator $\prod_{j=2}^ND_{\alpha_j}$ preserves real-rootedness (because $\alpha_j\ge 0$).
Hence if $p$ has only real zeros, the polynomial $\prod_{j=2}^ND_{\alpha_j}p$ has only real zeros as well. So, if $D p$ has a zero in $\zeta\not =0$ of multiplicity $l>1$, then Step 1 implies that  $\prod_{j=2}^ND_{\alpha_j}p$ has a zero at $\zeta$ of multiplicity $l+1$. Repeating the process, we deduce that $p$ has a zero at $\zeta$ of multiplicity $l+N$.
\end{proof}

\noindent
\textit{Step 3}. Let $(p_n)_n$ be the Brenke polynomials generated by $A$. If $p_n$ has a zero at $\zeta$ of multiplicity $l>1$, then $\zeta=0$.

\begin{proof}[Proof of Step 3]
For a real number $a$, write $a+\phi=(a+\phi_i)_{i=1}^N$ and denote $B_{a+\phi}$ for the hypergeometric function (\ref{h0f}) defined by the sequence $a+\phi$.

It is easy to check that $B'_\phi(z)=B_{1+\phi}/\theta_0$, where $\theta_0=\prod_{i=1}^N\phi_i\not =0$. Hence, if we consider the function $\tilde B_{1+\phi}$ defined by (\ref{tib}), we deduce that
$\tilde B_{1+\phi}=B_\phi$. Hence if we write
$(p_{n}^{1+\phi})_n$ for the Brenke polynomials generated by $A$ and associated to $B_{1+\phi}$, Remark \ref{urdc} gives that
\begin{equation}\label{smeb}
p_n=\frac{1}{\prod_{i=1}^N\phi_i}\Upsilon_{B_{1+\phi}}(p_{n}^{1+\phi}).
\end{equation}
The identity (\ref{oplh2}) gives
$$
\Upsilon_{B_{1+\phi}}=\prod_{i=1}^ND_{\phi_i}.
$$
Hence, since $\phi_i>0$, $i=1,\dots,N$, if $\zeta\not =0$ we deduce from Step 2 that  $p_{n}^{1+\phi}$ has at $\zeta$ a zero of multiplicity $l+N$. Repeating the process, we deduce that
$p_{n}^{n+\phi}$ has a zero at $\zeta$ of multiplicity $l+nN$. But this is impossible because the degree of $p_{n}^{n+\phi}$ is $n<l+nN$.
\end{proof}

\end{proof}

\medskip

\begin{proof}[Proof of Theorem \ref{pzgf0}]
Lemma \ref{jod} and Theorem \ref{tal} prove the simplicity property of the zeros of $p_n$.

The interlacing properties can be proved proceeding as in the proof of Lemma \ref{enz}, using again
the Obreshkov Theorem \ref{obre} (and the simplicity of the zeros).
\end{proof}

\bigskip

Actually, when $\phi_i>0$, $i=1,\dots, N$, and $A\in \rrp_{B_\phi}$ the zeros of the Brenke polynomials $(p_n)_n$ generated by $A$ and associated to $B_\phi$ seem to enjoy a lot of more properties. Here it is a trio of properties for which we have plenty of computational evidence but not a proof yet.

Since the zeros of $p_n$ are real, we denote by $\zeta^-_j=\zeta^-_j(\phi)$, $j=1,\dots, n_-$, for the negative zeros of $p_n$ ordered in decreasing size, where $n_-=n_-(\phi)$ denotes the number of negative zeros of $p_n$. Similarly, we denote by $\zeta^+_j=\zeta^+_j(\phi)$, $j=1,\dots, n_+$, for the positive zeros of $p_n$ ordered in increasing size, where $n_+=n_+(\phi)$ denotes the number of positive zeros of $p_n$.

\medskip

\noindent \textit{Invariance of $n_+$ and $n_-$ with respect to the parameter $\phi$}. We guess that the number of positive zeros of $p_n$ only depends on $A$ and nor on $B_\phi$. And so the same happens for the number of negative zeros.

\medskip

\noindent \textit{Monotonicity of the zeros of $p_n$ with respect to the parameters $\phi_i$}. We say that $\phi\preceq \psi$ if $\phi_i\le \psi_i$, for all $i=1,\dots ,N$.
Then, our conjecture is:
\begin{enumerate}
\item The $j$-th positive zero $\zeta_j^+$ is an increasing function of the parameter set $\phi$: if $\phi\preceq \psi$ then $\zeta^+_j(\phi)\le \zeta^+_j(\psi)$.
\item The $j$-th negative zero $\zeta_j^-$ is a decreasing function of the parameter set $\phi$: if $\phi\preceq \psi$ then $\zeta^-_j(\phi)\ge \zeta^-_j(\psi)$.
\end{enumerate}

\medskip

\noindent \textit{Interlacing properties of the zeros of $p_n$ for the parameters $\phi$ and $1+\phi$}. Our conjecture for the negative zeros is
$$
\zeta_{n_-}^-(1+\phi)<\zeta_{n_-}^-(\phi)<\zeta_{n_--1}^-(1+\phi)<\zeta_{n_--1}^-(\phi)<\dots <\zeta_{1}^-(1+\phi)<\zeta_{1}^-(\phi).\\
$$
And for the positive zeros
$$
\zeta_{1}^+(\phi)<\zeta_{1}^+(1+\phi)<\zeta_{2}^+(\phi)<\zeta_{2}^+(1+\phi)<\dots <\zeta_{n_+}^+(\phi)<\zeta_{n_+}^+(1+\phi).
$$

\subsection{The Appell-Dunkl case}
Appell-Dunkl polynomials are the particular case of Brenke polynomials defined by the Dunkl exponential, i.e.,
$B(z)=E_\mu(z)$, $\mu\not =-1,-2,\dots $ (see (\ref{emu})). A simple computation gives
\[
B(z)= \Ea(z) = \sum_{n=0}^\infty \frac{z^n}{\cc_{n,\mu}},
\]
with
\[
\cc_{n,\mu} =
\begin{cases}
  2^{2k}k!\,(\mu+1)_k, & \text{if $n=2k$},\\
  2^{2k+1}k!\,(\mu+1)_{k+1}, & \text{if $n=2k+1$}.
\end{cases}
\]

The (renormalized) Brenke polynomials generated by the formal power series $A$ and associated to $E_\mu(z)$ are defined by
\begin{equation}\label{adde}
A(z)E_\mu(xz)=\sum_{n=0}^\infty p_{n,\mu}(x)\frac{z^n}{\cc_{n,\mu}},
\end{equation}
and are called Appell-Dunkl polynomials because the operator $\Lambda_B$ (see (\ref{opgb})) is the Dunkl operator defined by
\begin{equation}
\label{cojo}
  \La f(x) = \frac{d}{dx}f(x) + \frac{2\mu+1}{2}
  \left(\frac{f(x)-f(-x)}{x}\right)
\end{equation}
(see, for instance, \cite{CDPV,Du,DPV,Ros}).
It is not difficult to check that the Appell-Dunkl polynomial $p_{n,\mu}$ is a monic polynomial of degree $n$ which satisfies
\begin{equation}
\label{coj}
\La p_{n,\mu}=(n+(\mu+1/2)(1-(-1)^n))p_{n-1,\mu}.
\end{equation}

The function $E_\mu(z)$, $\mu\not =-1,-2,\dots $, satisfies the hypothesis of  Theorem \ref{prin2} and hence, $\rrp_{E_\mu}\subset \rrp_{E_\mu}^p \subset \lp_0$.
But since $E_\mu(z)\not \in \lp I$, $\mu\not=-1/2$ (let us remind that $E_{-1/2}(z)=e^z$), we deduce from Theorem \ref{bpol} that $\rrp_{E_\mu}\subsetneq \lp_0$, $\mu\not\in \{-1/2,-1,-2,\dots \}$. Actually, it can be proved more: $\rrp_{E_\mu}^p\subsetneq \lp_0$. We just sketch the proof. Consider $A(z)=(z+1)^3\in \lp_0$ and the Brenke polynomials $(p_n)_n$ generated by $A$. Using (\ref{bps}), we have $p_n(x)=r_n(x)x^{n-3}$, $n\ge 3$, where $r_n$ is the following polynomial of degree $3$
$$
r_n(x)=x^3+3\frac{\cc_{n,\mu}}{\cc_{n-1,\mu}}x^2+3\frac{\cc_{n,\mu}}{\cc_{n-2,\mu}}x+\frac{\cc_{n,\mu}}{\cc_{n-3,\mu}}.
$$
A careful computation gives that
$$
\Delta (r_n)=\begin{cases} -2^43^3n^2(\mu+n/2)(2n\mu^2+(2n+1)\mu+n/2),&\mbox{$n$ even},\\
-2^53^3(n-1)\left(\mu+\frac{n+1}2\right)^2\left(2\mu(\mu+1)(2\mu+n+1)+\frac{n-1}2\right),&\mbox{$n$ odd}.
\end{cases}
$$
where $\Delta$ denotes the discriminant of the polynomial $r_n$ (see, for instance, \cite{vdw}). Hence
$$
\lim_{n\to \infty}\frac{\Delta (r_n)}{n^4}=-2^23^3(2\mu+1)^2.
$$
This gives that for $\mu\not=-1/2$ and $n$ big enough (depending on $\mu$), $\Delta(r_n)<0$, and then $p_n$ has two non-real zeros.
So, $A\not \in \rrp_{E_\mu}^p$ for $\mu\not\in \{-1/2,-1,-2,\dots \}$.

In general, the set $\rrp_{E_\mu}$ seems to strongly depend on $\mu$. However, using the next Lemma it is easy to describe the even functions in $\rrp_{E_\mu}$, $\mu>-1$.

\begin{Lem}\label{thzp} Assume that $A\in \lp_0 I$ with positive zeros, $b_n>0$, $n\ge 0$, and
consider the Brenke polynomials $(p_n)_n$ generated by $A$ and associated to $B$.
If $A\in \rrp_B$, then all the zeros of $p_n$ are positive for all $n\ge 0$.
\end{Lem}

\begin{proof}
Since the zeros of $A$ are positive, Theorem \ref{smac} implies that the sequence $(a_n)_n$ has alternating sign. Then, (\ref{bps}) shows that the coefficients of $p_n$ alternate sign as well. Since $p_n$ has only real zeros, they have to be positive.
\end{proof}

\begin{Cor}\label{TYY} Let $\mu$ and $A$ be a real number $\mu>-1$ and an even formal power series with $a_0=1$, respectively. Then $A\in \rrp_{E_\mu}$ if and only if $A\in\lp_0$.
In which case, if we write $(p_n)_n$ for the Appell-Dunkl polynomials generated by $A$, we have
\begin{enumerate}
\item If $A$ is a polynomial of degree $2k$ then for $0\le n\le 2k$, $p_n$ has simple zeros and the zeros of $p_{n-1}$ strictly interlace with the zeros of $p_{n}$. For $n\ge 2k$, we have
\begin{equation}\label{pesf}
p_n(x)=x^{n-2k}r_{n}(x)
\end{equation}
where $r_n$ is a polynomial of degree $2k$ with $r_n(0)\not=0$ and simple zeros, and the zeros of $p_n$ strictly interlace  the zeros of $p_{n+1}$ except for the common zero $x=0$.
\item If $A$ is not a polynomial then $p_n$ has simple zeros and the zeros of $p_n$ strictly interlace the zeros of $p_{n+1}$.
\end{enumerate}
Moreover, in that case $(1+az)A(z)\in\rrp_{E_\mu}$ for all $a\in\RR$.
\end{Cor}

\begin{proof}
Write $(p_n)_n$ for the Appell-Dunkl polynomials generated by $A$ (\ref{adde}). Since $A$ is even using (\ref{emu}), it is easy to see that
\begin{equation}\label{pqn}
p_{2n}(x)=\gamma_{2n,\mu}q_{n,\mu}((x/2)^2),\quad p_{2n+1}(x)=\frac{\gamma_{2n+1,\mu}x}{2(\mu+1)}q_{n,\mu+1}((x/2)^2),
\end{equation}
where $q_{n,\mu}$ are the Brenke polynomials generated by $A(\sqrt z)$ associated to $B_{1+\mu}$ (because $\I_\mu(z)=B_{1+\mu}((z/2)^2)$).

On the one hand, we have already prove that $\rrp_{E_{\mu}}\subset \lp$.

On the other hand, if $A\in \lp$ is even and $A(0)=1$, we have
$$
A(z)=e^{-az^2}\prod_{j=1}^\infty (1-\zeta_j^2z^2),
$$
with $a\ge 0$, $\zeta_j\in \RR$ and $\sum_j\zeta^2_j<+\infty$. Hence $A(\sqrt z)\in \lp I$ and has positive zeros.
Lemma \ref{thzp} gives that all the zeros of $q_{n,\mu}$ and $q_{n,\mu+1}$ are positive. (\ref{pqn}) then proves that the zeros of $p_n$ are real, and hence $A\in \rrp_{E_{\mu}}$.

Because of Theorem \ref{pzgf0}, we have only to prove the interlacing properties of the zeros.
Since $A$ is even, we have that $p_{2n}$ are even polynomials and $p_{2n+1}$ are odd polynomials. Hence, using (\ref{cojo}) and (\ref{coj}), we have
\begin{equation}\label{for}
p_{n-1}(x)=\begin{cases} \displaystyle \frac{1}{n}p_{n}'(x),&\mbox{if $n$ is even},\\
\displaystyle \frac{1}{n+2\mu+1}\left(p_{n}'(x)+\frac{2\mu+1}{x}p_{n}(x)\right),&\mbox{if $n$ is odd}.
\end{cases}
\end{equation}
Hence, the interlacing properties of the zeros of $p_n$ and $p_{n-1}$ is straightforward if $n$ is even.

If $n$ is odd, we proceed as follows. From Lemma \ref{jod}, we know that $x=0$ is the only zero of $p_n$ which might have multiplicity bigger than $1$. (\ref{for}) shows that if $p_n$ has a zero at $x=0$ of multiplicity $l$ then $p_{n-1}$ has also a zero of multiplicity $l-1$. Take now two consecutive zeros $\xi_1<\xi_2$ of $p_n$, $\xi_1,\xi_2\not =0$ and of the same sign. (\ref{for}) shows that $p_{n-1}(\xi_1)p_{n-1}(\xi_2)<0$. This shows that the zeros of $p_{n-1}$ interlace the zeros of $p_n$.

Finally, the Appell-Dunkl polynomials generated by $(1+az)A(z)$ are $p_n(x)+\frac{a\gamma_{n,\mu}}{\gamma_{n-1,\mu}}p_{n-1}(x)$. Since $p_n$ and $p_{n-1}$ interlace their zeros, the Obreshkov Theorem \ref{obre} gives that all the zeros of $p_n(x)+\frac{a\gamma_{n,\mu}}{\gamma_{n-1,\mu}}p_{n-1}(x)$ are real.
\end{proof}

\medskip

\section{Asymptotic for Brenke polynomials II}\label{ult}
In this section, we find some more asymptotics for Brenke polynomials  which provide new equivalencies for the Riemann Hypothesis in terms of real-rootedness of related sequences of Brenke polynomials.

For a formal power series $C(z)=\sum_{n=0}^\infty c_nz^n$, $c_0=1$ and $c_n\not =0$, $n\ge 0$, we can extend the linear operator $\Lambda _C$ (\ref{opgb}) from the linear space of polynomials to that of formal power series
as follows:
$$
\Lambda_C\left(\sum_{n=0}^\infty d_nz^n\right)=\sum_{n=0}^\infty d_{n+1}\frac{c_n}{c_{n+1}}z^n.
$$

\begin{Theo}\label{sasi}
Let $A$, $B$ and $C$ be normalized formal power series satisfying:
\begin{equation}\label{cabc}
a_n,c_n\not =0, \mbox{for all $n\ge 0$}, \quad \lim_n\frac{a_{n-1}a_{n+1}}{a_n^2}=1,\quad\lim_n\frac{c_{n-1}c_{n+1}}{c_n^2}=1.
\end{equation}
Write $(p_{n,s})_n$ for the Brenke polynomials generated by the formal power series $\frac{c_s}{a_s}\Lambda_C^sA$ associated to $B$. Then
\begin{equation}\label{asa}
\lim_s\frac{a_sc_{n+s}}{a_{n+s}c_s}p_{n,s}\left(\frac{a_{n+s+1}c_{n+s}}{a_{n+s}c_{n+s+1}}z\right)=r_n(z),
\end{equation}
where $(r_n)_n$ are the Brenke polynomials generated by $C$ and associated to $B$ (\ref{bps}).
Moreover, if $B\in \lp I$ and $C\in \lp$,  for every $n\ge 0$, there exists $s_n\ge 0$ such that
the polynomial $p_{n,s}$ has only real zeros for $s\ge s_n$.
\end{Theo}

\begin{proof}
If we apply $s$ times the operator $\Lambda_C$ to the formal power series $A$, we get (after normalization)
$$
\frac{c_s}{a_s}\Lambda_C^sA(z)=\frac{c_s}{a_s}\sum_{n=0}^\infty a_{n+s}\frac{c_n}{c_{n+s}}z^n.
$$
If we set
$$
\hat a_n=\frac{a_n}{c_n},
$$
according to (\ref{bps}), we get
\begin{align*}
p_{n,s}(z)&=\frac{c_s}{a_s}\sum_{j=0}^na_{n-j+s}\frac{c_{n-j}}{c_{n-j+s}}b_jz^j\\
&=\frac{1}{\hat a_s}\sum_{j=0}^n \hat a_{n-j+s}c_{n-j}b_jz^j.
\end{align*}
From (\ref{cabc}), it follows easily that the sequence $(\hat a_n)$ also satisfies
$$
\hat a_n\not =0, n\ge0,\quad \lim_n\frac{\hat a_{n-1}\hat a_{n+1}}{\hat a_n^2}=1.
$$
Hence, according to Lemma \ref{lems}, by writing $\mu_n=a_nc_{n+1}/(a_{n+1}c_n)$ we have
\begin{equation}\label{asax}
\lim_n\frac{\hat a_{n-j}}{\hat a_n\mu_n^j}=1,\quad \mbox{for all $j\ge 0$}.
\end{equation}
Writing
$$
\frac{\hat a_s}{\hat a_{n+s}}p_{n,s}(z/\mu_{n+s})=\sum_{j=0}^n \frac{\hat a_{n-j+s}}{\hat a_{n+s}\mu _{n+s}^j}c_{n-j}b_jz^j,
$$
and using (\ref{asax}), we get
$$
\lim_s\frac{\hat a_s}{\hat a_{n+s}}p_{n,s}(z/\mu_{n+s})=\sum_{j=0}^n c_{n-j}b_jz^j,
$$
which it is the asymptotic (\ref{asa}), since the polynomials in the right hand side of the previous formula are the Brenke polynomials generated by $C$ and associated to $B$.

As a consequence, if $B\in \lp I$ and $C\in \lp$, we deduce using Theorem \ref{bpol} that the polynomials $r_n$ have only real zeros. And so for every $n\ge 0$, there exists $s_n\ge 0$ such that
the polynomial $p_{n,s}$ has also only real zeros for $s\ge s_n$.
\end{proof}

\medskip

Take now $B(z)=e^z$ and $C(z)=\pFq{0}{1}{-}{\alpha+1}{-z}$, $\alpha>-1$, so that (see Remark~\ref{urdc})
$$
\Lambda_C=\frac{d}{dz}\left(\alpha+z\frac{d}{dz}\right)= (1+\alpha)\frac{d}{dz}+z\frac{d^2}{dz^2}.
$$
The Brenke polynomials $p_{n,s}$ generated by $\frac{c_s}{a_s}\Lambda_C^sA$ associated to $B$ are
$$
p_{n,s}(z)=\frac{(-1)^s}{(\alpha+1)_ss!\,a_s}\sum_{j=0}^n\frac{(n-j+1)_s(\alpha+n-j)_s}{j!\,(n-j)!} a_{n-j+s}z^j.
$$
And so, according to the asymptotic (\ref{asa}) (we assume that $A$ satisfies the hypothesis of Theorem \ref{sasi})
\begin{align}\label{ass}
\lim_s&\frac{(-1)^na_s}{a_{n+s}(\alpha+s+1)_n(s+1)_n}p_{n,s}\left(-(\alpha+n+s+1)(n+1+s)\frac{a_{n+s+1}}{a_{n+s}}z\right)\\\nonumber
&=\sum_{j=0}^n\frac{(-1)^j}{j!\,(n-j)!\,(\alpha+1)_j}z^{n-j}=\frac{z^n}{(\alpha+1)_n}L_n^{\alpha}(1/z),
\end{align}
where $L_n^\alpha$ is the $n$-th Laguerre polynomial.

Particularizing for $A=\varsigma$ (\ref{var}) we have Corollary \ref{esf1} (in the Introduction).

\begin{proof}[Proof of Corollary \ref{esf1}]

(2) $\Rightarrow$ (1).
For $s=0$, the polynomials $(p_{n,0}^\alpha)_n$ are the Jensen polynomials for $\varsigma$ (\ref{jenp}) and hence, RH is equivalent to $p_{n,0}^\alpha$ having only real zeros for all $n\ge 0$.

(1) $\Rightarrow$ (2).
We have mentioned in the Introduction that RH is equivalent to $\varsigma\in \lp I$. It follows from Theorem \ref{hedc} that the operator $\Lambda_C$ is stable.
Since $\varsigma\in \lp$ (we are assuming RH), we have from Corollary \ref{sccd} that $\Lambda_C^s\varsigma\in \lp$. Since $(p_{n,s}^\alpha)_n$ are the Brenke polynomials generated by $\Lambda_C^s\varsigma\in \lp$ associated to $B(z)=e^z\in \lp I$, we deduce from Theorem \ref{bpol} that $\Lambda_C^s\varsigma\in \rrp_B=\lp$ and hence the polynomials $p_{n,s}^\alpha$ have only real zeros for all $n\ge 0$.

The second part of the corollary is an easy consequence of the asymptotic (\ref{ass}).
\end{proof}

\bigskip

Proceeding in a similar form, we can prove the following dual result for Theorem~\ref{sasi}.
\begin{Theo}\label{sasi2}
Let $A$, $B$ and $C$ be normalized formal power series satisfying
\begin{equation}\label{cabc3}
b_n,c_n\not =0, \mbox{for all $n\ge 0$}, \quad \lim_n\frac{b_{n-1}b_{n+1}}{b_n^2}=1,\quad\lim_n\frac{c_{n-1}c_{n+1}}{c_n^2}=1.
\end{equation}
Write $(q_{n,s})_n$ for the Brenke polynomials generated by $A$ and associated to $\frac{c_s}{b_s}\Lambda_C^sB$. Then
\begin{equation}\label{asajj}
\lim_s\frac{b_sc_{n+s}}{b_{n+s}c_s}\left(\frac{b_{n+s+1}c_{n+s}}{b_{n+s}c_{n+s+1}}z\right)^nq_{n,s}\left(\frac{b_{n+s}c_{n+s+1}}{b_{n+s+1}c_{n+s}}\frac{1}{z}\right)=q_n(z),
\end{equation}
where $(q_n)_n$ in the right hand side of the previous identity are the Brenke polynomials generated by $C$ and associated to $A$.
Moreover, if $A\in \lp I$ and $C\in \lp$, the polynomials $q_n$ have only real zeros. And so for every $n\ge 0$, there exists $s_n\ge 0$ such that
the polynomial $q_{n,s}$ has also only real zeros for $s\ge s_n$.
\end{Theo}

When $C(z)=e^z$, we have $\Lambda_C=d/dz$, and hence the previous asymptotic leads to Corollary \ref{ggorz}, which provides an alternative proof of Theorem 1 in \cite{GORZ}.

\bigskip

Particularizing for $A(z)=\pFq{0}{1}{-}{\alpha+1}{-z}$ and $B(z)=\varsigma(z)$ we have Corollary \ref{esf2} (in the Introduction).

\begin{proof}[Proof of Corollary \ref{esf2}]

(2) $\Rightarrow$ (1).
For $s=0$, the polynomials $(q_{n,0}^\alpha)_n$ are the Brenke polynomials generated by $A(z)=\pFq{0}{1}{-}{\alpha+1}{-z}$ and associated to $B(z)=\varsigma(z)$. Hence $A\in \rrp_\varsigma$. Lemma~\ref{sdc} and Corollary~\ref{ejh} give that $\varsigma\in \rrp_A=\lp$. And $\varsigma\in\lp$ is equivalent to the Riemann hypothesis.

(1) $\Rightarrow$ (2).
Taking $s\ge 0$ and assuming that RH is true, we have that $\varsigma^{(s)}\in \lp I$. Since $B(z)=\varsigma^{(s)}(z)/\varsigma^{(s)}(0)$ satisfies (\ref{supl}), Theorem \ref{prin2} implies $\lp=\rrp_B$.
Since $A\in \lp=\rrp_B$, we can conclude that the zeros of $q_{n,s}^\alpha$ have to be real.

The second part of the corollary is an easy consequence of the asymptotic (\ref{asajj}).
\end{proof}

\bigskip

\noindent
\textit{Acknowledges.} I want to thank J. Arias de Reyna (my \textit{maestro}) and M. P\'erez and J.L. Varona (my friends and colleagues) for the discussions during the preparation of this paper, and for their careful reading of earlier versions of this paper.



\end{document}